\documentclass[leqno]{article}
\usepackage{amsmath,amsfonts,amsthm,amssymb,indentfirst}

\newcommand{\<}{\kern.0833em}

\newtheorem{theorem}{Theorem}[section]
\newtheorem{lemma}[theorem]{Lemma}
\newtheorem{corollary}[theorem]{Corollary}
\newtheorem{proposition}[theorem]{Proposition}
\newtheorem{definition}[theorem]{Definition}
\newtheorem{question}[theorem]{Question}

\newcommand{\fb}{\mathbf}
\newcommand{\R}{\mathbb R}
\newcommand{\Z}{\mathbb Z}
\newcommand{\Q}{\mathbb Q}
\newcommand{\Sym}[1][(\Omega)]{\mathrm{Sym}#1}
\newcommand{\Se}[1][(\Omega)]{\mathrm{Self}#1}
\newcommand{\En}[1][(V)]{\mathrm{End}_K #1}
\newcommand{\Eq}[1][(\Omega)]{\mathrm{Equiv} #1}
\newcommand{\Rel}[1][(\Omega)]{\mathrm{Rel}#1}
\newcommand{\cd}[1][\Omega]{{\mathrm{card}(#1)}}
\newcommand{\V}{\fb{V}}
\newcommand{\Pw}[1][\Omega]{{\fb{P}(#1)}}
\newcommand{\Pf}[1][\Omega]{{\fb{P}_\mathrm{fin}(#1)}}
\newcommand{\cP}{{\scriptstyle\>\coprod\>}}
\newcommand{\op}{^{\mathrm{op}}}
\newcommand{\strt}[1][1.7]{\vrule width0pt height0pt depth#1pt}% ^ strut

\makeatletter
\newcommand{\xlabel}{\stepcounter{equation}
  \gdef\@currentlabel{\p@equation\theequation}{\rm(\@currentlabel)}}
\makeatother
\newenvironment{xlist}
  {\begin{list}{\xlabel}
    {\setlength{\rightmargin}{20pt}
     \setlength{\leftmargin}{37pt}
     \setlength{\labelsep}{20pt}
     \setlength{\labelwidth}{20pt}}}
  {\end{list}}

\makeatletter
\let\@zzmaketitle\@maketitle
\def\@maketitle{%
\let\@zznull\null
\setbox\@tempboxa\hbox to\hsize{\hss\raise 20mm\hbox{\texttt{Comments, corrections, and related references welcomed, as always!}}\hss}%
\ht\@tempboxa=0pt\dp\@tempboxa=0pt%
\def\null{\box\@tempboxa}%
\@zzmaketitle
\let\null\@zznull
}
\makeatother

\begin{document}

\title{Some results on embeddings of
algebras, after~de~Bruijn~and~McKenzie%
\thanks{
%  	 *P,cP...
Secondary: 06Bxx, 08B20, 16S50, 18A99, 20B07, 20M20, 54Hxx.
% lats:xx=05,15-0 free   end(M) functrs   Sym    Se  top-alg:xx=11-5
This preprint is readable online at
http://math.berkeley.edu/\protect\linebreak[0]%
{$\!\sim$}gbergman\protect\linebreak[0]%
/papers/\,, and arXiv:math/0606407\,.
The former version is likely to be updated more frequently than
the latter.
}}
\author{George M. Bergman}
\maketitle

\begin{abstract}
In 1957, N.\,G.\,de~Bruijn showed that
the symmetric group $\Sym$ on an infinite set $\Omega$
contains a free subgroup on $2^{\cd}$ generators, and
proved a more general statement, a sample consequence
of which is that for any group $A$ of cardinality $\leq\cd,$
the group $\Sym$ contains a coproduct of $2^{\cd}$ copies of $A,$ not
only in the variety of all groups, but in
any variety of groups to which $A$ belongs.
His key lemma is here generalized to an arbitrary variety
of algebras $\V,$ and formulated as a statement about
functors $\fb{Set}\rightarrow\V.$
From this one easily obtains analogs of the results stated above
with ``group'' and $\Sym$
replaced by ``monoid'' and the monoid $\Se$ of endomaps of $\Omega,$
by ``associative $\!K\!$-algebra'' and the $\!K\!$-algebra
$\En$ of endomorphisms of a $\!K\!$-vector-space $V$ with basis
$\Omega,$ and by ``lattice''
and the lattice $\Eq$ of equivalence relations on $\Omega\<.$
It is also shown, extending another result from de~Bruijn's 1957
paper, that each of $\Sym,$ $\Se$ and $\En$ contains
a coproduct of $2^{\cd}$ copies of itself.

That paper also gave an example of a group of cardinality
$2^{\cd}$ that was {\em not} embeddable in $\Sym,$ and R.\,McKenzie
subsequently established a large class of such examples.
Those results are shown here to be instances of a
general property of the lattice of solution sets in $\Sym$ of
sets of equations with constants in $\Sym.$
Again, similar results -- this time of varying strengths -- are obtained
for $\Se,$ $\En,$ and $\Eq,$ and also for the monoid $\Rel$
of binary relations on~$\Omega\<.$

Many open questions and areas for further investigation are noted.
\end{abstract}
% - - - - - - - - - - - - - - - - - - - - - - - - - - - - - -

\section{Conventions, and outline.}\label{S.intro}

Throughout this note, $\Omega$ will be an infinite set.
Each ordinal (in particular, each
natural number) is understood to be the set of all smaller
ordinals; the set of all natural numbers is denoted $\omega.$
Functions, including elements of permutation
groups, will be written to the left of their
arguments and composed accordingly.
The word ``algebra'' will be used in the sense of
general algebra (universal algebra), except in the combination
``$\!K\!$-algebra'', which will always mean an associative unital
algebra in the sense of ring theory, over a field $K$ assumed
fixed throughout this note.
In those contexts, $V$ will denote a vector space
with basis $\Omega$ over that field $K.$

In \S\S\ref{S.free}-\ref{S.cP+F} we develop results to the effect
that algebras arising as values of certain sorts of functors
can be embedded in certain infinite direct product algebras, and
obtain, as immediate corollaries, results on embeddability of
groups, monoids, $\!K\!$-algebras, and lattices
in the group $\Sym,$ the monoid $\Se,$ the
$\!K\!$-algebra $\En,$ and the lattice $\Eq$ respectively
(all defined as in the abstract).
The remaining sections obtain results specific to embeddings in one or
another of those four structures, and in the monoid $\Rel.$
In \S\ref{S.jiggle} (and two appendices,
\S\S\ref{S.jiggleEn}-\ref{S.paths}) it is shown that one
can embed into each of the first
three of these algebras a coproduct of $2^{\cd}$
copies of that same algebra, while \S\S\ref{S.notSym}-\ref{S.notEn}
obtain restrictions on algebras $A$ embeddable in these five algebras,
in terms of order-properties of chains of solution sets of
systems of equations in $A.$
\S\ref{S.others} suggests some ways in which the
results of this note might be extended.

For some further unusual properties of $\Sym$ and of some of the other
structures here considered, cf.~\cite{Sym_Omega:1}, \cite[\S6]{P_vs_cP},
\cite{Sym_Omega:2}, and works referred to in those papers.

I am indebted to
Anatole Khelif
for inadvertently bringing these questions to my attention, to
Zachary Mesyan
for a careful reading of the first draft, to
Andreas Blass and Vladimir Tolstykh
for helpful information on the literature,
and to the referee for several useful suggestions and corrections.

\section{Free algebras.}\label{S.free}

Recall that by Cayley's Theorem, every group of cardinality
$\leq\cd$ can be embedded in the symmetric group $\Sym$ on $\Omega\<;$
in particular, $\Sym$ contains free groups of all ranks $\leq\cd.$
An obvious question is whether it contains larger free groups,
for example, a free group of rank $\cd[\Sym]=2^{\cd}.$
In \cite{debruijn}, de~Bruijn answered this question affirmatively
by a method which also gave embeddings of many
interesting nonfree groups in $\Sym.$
We will begin by illustrating his key trick, concerning subalgebras
of direct products, in the case of free algebras, making the trivial
generalization from the variety of groups to arbitrary varieties
of algebras.
In the next section, his more general statement will be motivated,
reformulated in functorial terms, and generalized still further.

In these two sections, $\V$ will be any variety of finitary algebras.
(``Finitary'' means that every operation has finite arity,
but does not exclude varieties with infinitely many operations;
for example, modules over an infinite ring.)

\begin{proposition}[{\rm cf.~\cite{debruijn}}]\label{P.free}
Let $H$ be the free algebra on $\aleph_0$ generators in a
variety $\V.$
Then the direct product
algebra $H^{\cd}$ has a subalgebra free on $2^{\cd}$ generators.
\end{proposition}\begin{proof}
Let $\Pw$ denote the power set of $\Omega,$ and $\Pf\subseteq\Pw$ the
set of finite subsets of $\Omega\<.$
Then $\cd[\Pw]=2^{\cd}$ and $\cd[\Pf]=\cd,$ so
it will suffice to find a $\!\Pw\!$-tuple of elements of $H^{\Pf}$
that satisfies no relations other than the identities of $\V.$

For each $s\in\Pf,$ let $H_s$ denote the factor indexed
by $s$ in our product $H^{\Pf},$ and
let us pick $2^{\cd[s]}$ of the $\aleph_0$ free generators of $H_s,$
denoting these $x_{s,t},$ with $t$ ranging over the subsets of $s.$
For every $r\in\Pw,$ let $X_r$ be the element of $H^{\Pf}$
which, for each $s\in\Pf,$ has $\!H_s\!$-component $x_{s,r\cap s}.$

We claim that the $\!\Pw\!$-tuple of elements $(X_r)_{r\in\Pw}$
satisfies no relations other than identities of $\V.$
Indeed, since $\V$ is finitary, any relation satisfied by
these elements involves only finitely many of them;
let $R(X_{r_1},\dots,X_{r_n})$ be such a relation,
where $r_1,\dots,r_n$ are distinct elements of $\Pw.$
Choose a finite subset $s\subseteq\Omega$
such that $r_1\cap s,\dots,r_n\cap s$ are distinct.
Then the $s\!$-components of $X_{r_1},\dots,X_{r_n},$
namely $x_{s,r_1\cap s},\dots,x_{s,r_n\cap s},$ are independent
indeterminates in $H_s,$ so projecting the
relation $R(X_{r_1},\dots,X_{r_n})$ that we assumed to hold in
$H^{\Pf}$ onto the component $H_s$ of that product,
we see that it is an identity of~$\V.$
\end{proof}

\begin{corollary}\label{C.free}
\textup{(i)}~ The symmetric group $\Sym$ on $\Omega$ has subgroups
free on $2^{\cd}$ generators in every variety $\V$ of groups.\\
\textup{(ii)}~ The monoid $\Se$ of endomaps of $\Omega$ has submonoids
free on $2^{\cd}$ generators in every variety $\V$ of monoids.\\
\textup{(iii)}~ The endomorphism ring $\En$ of the $\!K\!$-vector
space $V$ with basis $\Omega$ has $\!K\!$-subalgebras free on
$2^{\cd}$ generators in every variety $\V$ of
associative $\!K\!$-algebras.\\
{\rm(iv)}~ The lattice $\Eq$ of equivalence
relations on $\Omega$ has sublattices free on
$2^{\cd}$ generators in every variety $\V$ of lattices.
\end{corollary}\begin{proof}
Cayley's Theorem shows that every group of
cardinality $\leq\cd$ is embeddable in $\Sym,$
and similar arguments give embeddings
of all monoids with $\leq\cd$ elements in $\Se$
and of all associative \mbox{$\!K\!$-algebras} of vector-space dimension
$\leq\cd$ in $\En.$
The corresponding statement for embeddability of lattices
in $\Eq$ is Whitman's Theorem \cite{whitman}.
(Whitman does not explicitly say there that if the
given lattice $L$ is infinite, then the set on which he represents it
has cardinality $\leq\cd[L],$ but this
can be verified from his construction;
or one can deduce the possibility of an embedding with this
cardinality condition from the embeddability result without it.)
In particular, each of these structures contains a copy of the free
algebra $H$ on $\aleph_0$ generators in any subvariety $\V$ of
the given variety, since that free algebra is countable, or in the
$\!K\!$-algebra case, countable dimensional.

Moreover, each of the four algebras named contains a
$\!\cd\!$-fold direct product of copies of itself.
To see this, let us write the given set $\Omega$ as a disjoint union of
$\cd$ subsets of cardinality $\cd,$ $\Omega=\bigcup_{i\in\cd}\Omega_i.$
Then within $\Sym,$ the subgroup of consisting those permutations that
respect each $\Omega_i$ is such a direct product, and
the analogous statement holds in the monoid $\Se.$
In the $\!K\!$-algebra case we similarly use the algebra
of endomorphisms that carry the span of each subset $\Omega_i$ of our
basis $\Omega$ of $V$ into itself,
and in the lattice case, the sublattice
of equivalence relations that relate members of each $\Omega_i$
only with other members of $\Omega_i.$

Since we saw in the first paragraph that each of our objects
contains a free algebra $H$ on countably many generators in
the (arbitrary) subvariety $\V,$ it follows that it
will contain a copy of the product algebra $H^{\cd},$ and
Proposition~\ref{P.free} now gives the desired conclusions.
\end{proof}

Note that the embedding of $\Eq^{\cd}$ in $\Eq$ used
in the second paragraph of the above proof takes the least element
of $\Eq^{\cd}$ to the least element of $\Eq,$ but does not
take the greatest element to the greatest element.
(It takes that element to the relation whose equivalence classes are
the sets $\Omega_i.)$
With a little more work, however, one can get the an embedding
that respects both greatest and least elements.

Namely, fix any $p\in\Omega,$ and let
$\Omega'\subseteq\Omega^{\cd}$ be the set of elements $(x_i)_{i\in\cd}$
such that $x_i=p$ for all but finitely many $i.$
We see that $\cd[\Omega']=\cd,$ so it will suffice to
embed $\Eq^{\cd}$ in $\Eq[(\Omega')].$
We do this by taking each $\!\cd\!$-tuple $(\alpha_i)_{i\in\cd}$
$(\alpha_i\in\Eq)$ to the relation $\alpha\in\Eq[(\Omega')]$
such that $((x_i),(y_i))\in\alpha$ if and only if
$(x_i,y_i)\in\alpha_i$ for all $i.$
It is straightforward to show that this is a lattice embedding
which indeed respects least and greatest elements.

So Corollary~\ref{C.free}(iv) also holds for lattices with
greatest and/or least element, and mappings respecting these elements.
For brevity, I will not mention lattices with this additional structure
in subsequent sections, except when points come up
where I notice that what we can prove about lattices with such
structure differs from what we can prove for lattices without it.
(Incidentally, both the above embedding and the one in the proof of
Corollary~\ref{C.free}(iv) also respect {\em infinitary}
meets and joins; but this is not relevant to our embedding results,
since those require that the algebra operations used be finitary.)

\section{Coproducts and functors.}\label{S.cP+F}

The free algebra on a set $\Omega$ in a variety $\V$
is the {\em coproduct} in $\V$ of an $\!\Omega\!$-tuple
of copies of the free algebra on one generator.
To start the ball of generalization rolling, let us note how to
extend the proof of Proposition~\ref{P.free} to the case where free
algebras are replaced by coproducts of copies of an arbitrary algebra.

In our proof of Proposition~\ref{P.free}, we chose in each copy
$H_s$ of $H$ a $\!\Pw[s]\!$-tuple $(x_{s,t})_{t\subseteq s}$ of
distinct members of our $\!\aleph_0\!$-tuple of free generators.
This time, let $H$ be the coproduct in $\V$ of $\aleph_0$ copies of a
fixed algebra $A,$ and let us take for each $s\in\Pf$
a $\!\Pw[s]\!$-tuple $(p_{s,t})_{t\subseteq s}$ of distinct members
of the $\!\aleph_0\!$-tuple of coprojection maps $A\rightarrow H_s$
defining the coproduct structure.
We can then define, for each $r\in\Pw,$ a map $P_r:A\rightarrow H^{\Pf}$
by letting the composite of $P_r$ with each projection
$H^{\Pf}\rightarrow H_s$ be $p_{s,r\cap s}.$
With these adjustments, the proof of
Proposition~\ref{P.free} goes over, and we likewise get the
corollary that if $A$ is any group, monoid, associative
$\!K\!$-algebra, or lattice, of cardinality, respectively
$\!K\!$-dimension, $\leq\cd,$ and $\V$ any variety of groups,
monoids, \mbox{$\!K\!$-algebras} or lattices containing $A,$ then
the group $\Sym,$ the monoid $\Se,$ the
$\!K\!$-algebra $\En$ or the lattice $\Eq$
contains a coproduct in $\V$ of $2^{\cd}$ copies of $A.$\vspace{6pt}

To suggest the next level of generalization,
let me give a more or less random concrete example.
Let $A$ be the group presented by two generators, $x$ and $y,$
and the relations saying that the generator $x$ has exponent $2,$
and commutes with the element obtained by conjugating it by
the square of the generator $y:$
\begin{xlist}\item\label{x.ord2}
$x^2\ =\ 1,$
\end{xlist}
\begin{xlist}\item\label{x.xy2}
$x\,(y^2\,x\,y^{-2})\ =\ (y^2\,x\,y^{-2})\,x.$
\end{xlist}

Now if $I$ is any index-set, let $F(I)$ be the group presented
by generators $x_i,\ y_i$ $(i\in I)$ subject to the relations
\begin{xlist}\item\label{x.ord2i}
$x_i^2\ =\ 1$\quad$(i\in I),$
\end{xlist}
\begin{xlist}\item\label{x.xyyijk}
$x_i\,((y_jy_k)\,x_i\,(y_jy_k)^{-1})\ =\ %
((y_jy_k)\,x_i\,(y_jy_k)^{-1})\,x_i$\quad$(i,j,k\in I,$ not necessarily
distinct).
\end{xlist}
Looking at the $i=j=k$ case of these relations,
we see that for each $i,$ there is a homomorphism
$A\rightarrow F(I)$ acting by $x\mapsto x_i,\ y\mapsto y_i.$
This is in fact an embedding, for we also
see from (\ref{x.ord2})-(\ref{x.xyyijk}) that there exists a
homomorphism $F(I)\rightarrow A$ mapping all $x_i$ to $x$
and all $y_i$ to $y,$ which gives
a left inverse to each of the preceding homomorphisms.
(Note, incidentally, that the choices made in~(\ref{x.ord2i})
and~(\ref{x.xyyijk}), to turn $x^2$ to $x_i^2,$
but $y^2$ to $y_j\,y_k,$ were somewhat arbitrary: other choices
would have led to these same conclusions, so the
relations~(\ref{x.ord2}) and~(\ref{x.xy2}) did not uniquely
determine~(\ref{x.ord2i}) and~(\ref{x.xyyijk}).)

Though the group $F(I)$ is generated by an $\!I\!$-tuple of embedded
homomorphic images of $A,$ our presentation does not make
it a coproduct of those subgroups, since
the relations~(\ref{x.xyyijk}) relate elements from
{\em different} copies of $A;$ nor does it make it their coproduct
in some subvariety of groups, since~(\ref{x.ord2i})
and~(\ref{x.xyyijk}) do not
describe identities satisfied by all elements of $F(I).$
We see, however, as for coproducts, that any map of index-sets
$I\rightarrow J$ induces a group homomorphism $F(I)\rightarrow F(J),$
making $F$ a functor from sets to groups.

We shall find below that the idea of Proposition~\ref{P.free}
can be used to show that $F(\aleph_0)^{\cd}$
contains a copy of $F(2^{\cd}),$ and that the
corresponding statement holds with the variety of groups replaced
by any variety $\V$ of finitary algebras, and~(\ref{x.ord2i})
and~(\ref{x.xyyijk}) by any such system of ``relations parametrized
by families of indices''.

We could give a careful formulation of this concept of
a ``parametrized system of generators and relations''.
Fortunately, we do not have to, for we shall see that the concept
is equivalent to one that can be defined in a simpler way.
We noted above that the construction $F$ was a functor
$\fb{Set}\rightarrow\V;$ and it clearly satisfies
\begin{xlist}\item\label{x.genA}
For every set $I,$ the algebra $F(I)$ is generated by the
union of the images of $A=F(1)$ under the homomorphisms induced
by all maps $1\rightarrow I.$
\end{xlist}

I claim, conversely, that any functor $\fb{Set}\rightarrow\V$
satisfying~(\ref{x.genA}) corresponds to a system of algebras
determined by ``generators and relations with parameters''
in the sense suggested by the above discussion.
Indeed, given $F,$ let us take for
generator-symbols (corresponding to the $x$ and $y$ in
our group-theoretic example) any generating set $X$ for $F(1).$
For every set $I,$ every $i\in I,$ and every $x\in X,$
let us write $x_i$ for the image of $x$
under the map $F(1)\rightarrow F(I)$ induced by the
map $1\rightarrow I$ taking $0$ to $i.$
Then~(\ref{x.genA}) shows that $F(I)$ is generated by
\begin{xlist}\item\label{x.xi}
$\{x_i\mid x\in X,\ i\in I\}.$
\end{xlist}
To get relations, let us, for each natural
number $n,$ choose a set of relations presenting $F(n)$
in terms of the generators $x_i$ $(x\in X,\ i\in n),$
and let us turn each of these into a ``system of relations with
parameters'' by replacing the subscripts $0,\dots,n{-}1\in n$
on the generators appearing in each relation with symbols
$i_0,\dots i_{n-1}$ ranging over a general index-set $I.$

We see from the functoriality of $F$ that for any $I,$
the generators~(\ref{x.xi}) of $F(I)$
satisfy all instances of the system of relations so obtained.
To see that no more relations are needed, note that any
relation satisfied in $F(I)$ by the elements~(\ref{x.xi})
can involve only finitely many of these elements, say those coming
from the image of $F(n)$ under some one-to-one map
$n\rightarrow I,$ for some $n\in\omega.$
If $I\neq\emptyset,$ we can take $n>0,$ so that we may choose a left
inverse $I\rightarrow n$ to this map, and applying $F$ to it, we see
that the corresponding relation indeed holds in $F(n),$
and so is a consequence of the system of relations we have chosen.
If $I=\emptyset=0,$ then $n$ will also equal $0,$ and such a
map likewise exists, yielding the same conclusion.
Thus, the indicated system of generators and relations indeed
determines $F(I)$ for all $I.$
(Equations satisfied by the empty set of generators
correspond to relations on the set of zeroary operations of $\V,$
which hold in all $F(I)$ including $F(0).$
However, for the arguments below, we only need the values of $F(I)$
for nonempty index-sets $I,$ so nothing is lost if the reader prefers
to consider $F$ a functor from the category of nonempty sets
to $\V,$ and so avoid dealing with the case $I=\emptyset.)$

De~Bruijn~\cite{debruijn} proves his embeddability results for
what he calls ``symmetrically generated groups''.
On examination, these turn out to be precisely the values $F(I)$ of
group-valued functors $F$ satisfying~(\ref{x.genA}).
However, rather than stopping here, we may ask whether, in addition
to allowing relations like~(\ref{x.xyyijk}) that depend on more than one
parameter, we could allow this in our {\em generators} as well.
For example, suppose we associate to each set $I$ the
group $F(I)$ with generators
\begin{xlist}\item\label{x.xij}
$x_{ij}$ $(i,j\in I),$
\end{xlist}
subject to relations
\begin{xlist}\item\label{x.xxijk}
$x_{ij}x_{jk}=x_{jk}x_{ij}$ $(i,j,k\in I).$
\end{xlist}
It is again clear that maps among index-sets induce homomorphisms
among these groups, giving a functor $F:\fb{Set}\rightarrow\fb{Group},$
and that the corresponding statement is true for systems of algebras
of any variety $\V$ presented by generators and relations
similarly parametrized by multiple subscripts.
The resulting functors will not in general satisfy~(\ref{x.genA}), but
assuming the string of subscripts on each
generator is finite, they will satisfy
\begin{xlist}\item\label{x.genFfin}
For every set $I,$ the algebra $F(I)$ is the union of the images of
the homomorphisms $F(a): F(n)\rightarrow F(I),$ where $n$ ranges
over $\omega,$ and $a$ over all set-maps $n\rightarrow I.$
\end{xlist}
Conversely, it is straightforward to show, as before, that
the values of any functor satisfying~(\ref{x.genFfin}) arise from this
sort of presentation-with-parameters.

We can now give our generalization of Proposition~\ref{P.free}.
As indicated in the second paragraph of~\S\ref{S.free}, $\V$ denotes an
arbitrary fixed variety of finitary algebras.

\begin{theorem}[{\rm cf.~\cite[Theorem~3.1]{debruijn}}]\label{T.funct}
Let $F$ be a functor $\fb{Set}\rightarrow\V$
satisfying~\textup{(\ref{x.genFfin})}.
Then $F(\aleph_0)^{\cd}$ has a subalgebra
isomorphic to $F(2^{\cd}).$
\end{theorem}\begin{proof}
As in the proof of Proposition~\ref{P.free}, it suffices to
construct an embedding $h:F(\Pw)\rightarrow F(\aleph_0)^{\Pf};$
and again, we may specify such an $h$ by giving its
composites with the projections of $F(\aleph_0)^{\Pf}$ onto the
factors $F(\aleph_0)$ corresponding to each $s\in\Pf.$
For each such $s,$
let $c_s:\Pw\rightarrow\Pw[s]$ be defined by $r\mapsto r\cap s,$ and
let us choose an embedding $e_s:\Pw[s]\rightarrow\aleph_0,$
and take the composite of $h$ with the $\!s\!$th projection to be
$F(e_s c_s): F(\Pw)\rightarrow F(\aleph_0).$

To show that $h$ is an embedding, consider any two elements
$u\neq v\in F(\Pw).$
We claim there exists a component of $F(\aleph_0)^{\Pf}$ at which
$h(u)$ and $h(v)$ have distinct coordinates.

Indeed, from~(\ref{x.genFfin}) we can see that the images in
$F(\Pw)$ of homomorphisms $F(a): F(n)\rightarrow F(\Pw)$ induced by
maps $a:n\rightarrow\Pw$ $(n\in\omega)$ form a directed system
of subalgebras with union $F(\Pw).$
Hence $u$ and $v$ will together lie in such an image; so
let $u=F(a)(u_0),\ v=F(a)(v_0)$ for some $a:n\rightarrow\Pw$ and
$u_0, v_0\in F(n),$ necessarily distinct; here we may assume $n>0.$

Now choose $s\in\Pf$ such that
$a(0)\cap s,\dots,a(n{-}1)\cap s$ are distinct.
Then the composite map
$e_s\,c_s\,a: n\rightarrow\Pw\rightarrow\Pw[s]\rightarrow\aleph_0$ is
one-to-one, hence it has a left inverse.
Hence so does $F(e_s\,c_s\,a): F(n)\rightarrow F(\aleph_0);$
hence that is also one-to-one.
In particular, the images of $u_0$ and $v_0$ under the latter
map, which are the $\!s\!$-coordinates of $h(u)$ and $h(v),$
are distinct, as required.
\end{proof}

As before, we immediately get the particular embeddability results:

\begin{theorem}[{\rm cf.~\cite[Theorem~3.1]{debruijn}}]\label{T.funct'}
Suppose $F$ is a functor from $\fb{Set}$ to {\rm(i)}~the category
of groups, respectively {\rm(ii)}~the category of monoids,
{\rm(iii)}~the category of associative algebras over a field $K,$ or
{\rm(iv)}~the category of lattices; and suppose that $F$
satisfies~\textup{(\ref{x.genFfin})}, and
has the property that the cardinality of $F(\aleph_0)$ in
case {\rm(i)}, {\rm(ii)} or {\rm(iv)}, or its $\!K\!$-dimension
in case {\rm(iii),} is $\leq\cd.$
\textup{(}For instance, starting with
an algebra $A$ of cardinality or $\!K\!$-dimension $\leq\cd,$
one might define $F$ to be the functor associating to
every set $I$ the $\!I\!$-fold coproduct of copies of $A$ in some
fixed variety or quasivariety containing $A.)$

Then $F(2^{\cd})$ is embeddable in
{\rm(i)}~$\Sym,$ {\rm(ii)}~$\Se,$
{\rm(iii)}~$\En,$ or {\rm(iv)}~$\Eq,$ respectively.\qed
\end{theorem}

The reader may have noticed when we first proved
Proposition~\ref{P.free} that we did not really need the factors in our
product to be free of rank $\aleph_0;$ free objects of finite nonzero
ranks would do, as long as there were at least $\cd$ such factors of
rank greater than or equal to each natural number $N;$ and, similarly,
that in the proof of Theorem~\ref{T.funct}, we could have used a
direct product (with enough repetitions) of objects $F(n)$ for $n$
finite, instead of a power of $F(\aleph_0).$
However, it is not hard to verify in each of these cases that
the product of such a family would contain an embedded
copy of $F(\aleph_0)^{\cd},$ reducing these situations to that
of Theorem~\ref{T.funct}.
Let us record here the observation from which this follows.

\begin{lemma}[{\rm cf.~\cite{debruijn}}]\label{L.get*w}
Let $F$ be a functor $\fb{Set}\rightarrow\V$
satisfying~\textup{(\ref{x.genFfin})}.
Then $\prod_{0<n<\omega} F(n)$ has a subalgebra isomorphic to
$F(\aleph_0).$
\end{lemma}\begin{proof}
For each $n>0,$ let $f_n:\aleph_0=\omega\rightarrow n$ be the map taking
each natural number $r$ to $\min(r,n{-}1),$ and define
$f:F(\aleph_0)\rightarrow\prod_{0<n<\omega} F(n)$ to have
$F(f_n)$ as its $\!n\!$th coordinate, for each $n.$
An argument of the sort used in the proof of Theorem~\ref{T.funct}
shows that any two distinct elements of $F(\aleph_0)$ have distinct
projections in some $F(n),$ so $f$ is an embedding.
\end{proof}

It would be interesting to look for results similar to those
of this section for functors on categories other than $\bf{Set}.$
I leave these investigations to others, but give below
one such result I have noticed, and a couple
of examples of how it can be applied.

\begin{theorem}\label{T.F(R)}
Let $\bf{T.ord}$ be the category whose objects are totally
ordered sets, and whose morphisms are isotone maps \textup{(}maps
satisfying $x\leq y\implies a(x)\leq a(y)),$ and let
every ordinal, and likewise the set $\mathbb{R}$ of real numbers,
be regarded as objects of $\bf{T.ord}$ via their standard orderings.
Suppose $F: \fb{T.ord}\rightarrow \V$ is a functor satisfying the
analog of~\textup{(\ref{x.genFfin})} with ``totally ordered set''
for ``set'', and ``isotone maps'' for ``set-maps''.
Then $F(\omega)^{\aleph_0}$ has a subalgebra isomorphic to $F(\R).$
\end{theorem}\begin{proof}
The set $\Pf[\Q]$ of finite sets of rational numbers is
countable, so it suffices to embed $F(\R)$ in $F(\omega)^{\Pf[\Q]}.$
Given $s\in\Pf[\Q]$ whose distinct elements are $q_1<\dots<q_n,$
let $a_s: \R\rightarrow\omega$ be the isotone map which sends
each $r\in\R$ to the greatest $i\in\{1,\dots,n\}$ such that
$q_i\leq r$ if such an $i$ exists, and sends all $r<q_1$ to~$0.$
Let $h:F(\R)\rightarrow F(\omega)^{\Pf[\Q]}$ be the map whose
composite
with the projection indexed by each $s\in\Pf[\Q]$ is $F(a_s).$

It is not hard to see that we will be able
to complete the proof as we did that of Theorem~\ref{T.funct} if
for every finite set of real numbers $r_0<\dots<r_n,$ we can
find an $s\in\Pf[\Q]$ and an isotone map $b:\omega\rightarrow\R$
such that the map $b\<a_s:\R\rightarrow\omega\rightarrow \R$ fixes
$r_0,\dots,r_n.$
To do this, choose $q_1,\dots,q_n\in\Q$ so
that $r_{i-1}<q_i\leq r_i$ $(i=1,\dots,n),$
let $s=\{q_1,\dots,q_n\},$ and let $b$
take $i$ to $r_i$ for $i=0,\dots,n,$ and be extended
in an arbitrary isotone manner to larger $i.$
Thus, $a_s$ takes $r_i$ $(i=0,\dots,n)$ to $i,$ which $b$ takes
back to $r_i,$ as required.
\end{proof}

For a functor $F$ as in the above theorem, the algebras $F(I)$ will have
presentations by systems of generators and relations indexed by
finite sequences of subscripts
from $I,$ where the indices occurring in each generator
or relation may be constrained by inequalities of the form $i\leq j.$
A simple example is the functor associating to each totally
ordered set $I$ the group presented by
generators $x_i$ and $y_i$ $(i\in I)$ subject to the relations
\begin{xlist}\item\label{x.xyijR}
$x_i\<y_j = y_j\<x_i$ for $i\leq j$ in $I.$
\end{xlist}
Another is the functor taking each $I$ to
the (commutative) monoid presented by
generators $x_i$ $(i\in I)$ and relations
\begin{xlist}\item\label{x.xiR}
$x_i\<x_j = x_i = x_j\<x_i$ for $i\leq j$ in $I.$
\end{xlist}
In each of these cases, the object $F(\omega)$
is countable, hence embeddable in $\Sym[(\aleph_0)],$ respectively
$\Se[(\aleph_0)];$ hence, combining Theorem~\ref{T.F(R)} with the
method of proof of
Theorem~\ref{T.funct'}, we see that $F(\R)$ is also embeddable in
$\Sym[(\aleph_0)],$ respectively $\Se[(\aleph_0)].$
Likewise, the group algebra $K\<F(\R)$ for $F$ determined
by~(\ref{x.xyijR}),
and the monoid algebra $K\<F(\R)$ for $F$ determined by~(\ref{x.xiR}),
are the values at $\R$ of $\!K\!$-algebra-valued functors
satisfying the analog of~(\ref{x.genFfin}), and so
are embeddable in $\En$ for $V$ countable-dimensional.
I do not see any way to obtain these results from
Theorems~\ref{T.funct} and~\ref{T.funct'} themselves.

Let us note in connection with the group-theoretic
construction~(\ref{x.xyijR}) (and for some later uses)
that if $X$ is any set, and $R$
any symmetric reflexive binary relation on $X,$ and we form
the group $G$ presented by the generating set $X$ and the relations
\begin{xlist}\item\label{x.relB}
$x\<x'\ =\ x'x$ $((x,x')\in R),$
\end{xlist}
then distinct subsets of $X$ generate distinct subgroups of $G$
(clear by looking at the abelianization of $G),$ and elements
$x,\,x'\in X$ commute in $G$ if and only if $(x,x')\in R.$
To see the latter
statement, consider any $x\in X,$ and let $G_0$ be the group presented
as above, but using the set $X-\{x\}$ and the restriction of $R$
to that set.
Then $G$ can be described as an HNN extension~\cite{L+S}
of $G_0,$ obtained by adjoining an additional generator $x,$ whose
conjugation action is specified on the subgroup generated by
$\{y\in X-\{x\}\mid (x,y)\in R\}$ as the identity map.
By the structure of HNN extensions, conjugation by $x$ fixes precisely
the elements of that subgroup, giving the asserted characterization of
the commuting pairs of elements of $X.$

\section{Each of $\Sym,$ $\Se$ and $\En$ contains a coproduct
of copies of itself.}\label{S.jiggle}

In the case of Theorem~\ref{T.funct'} where $F(I)$ is the
$\!I\!$-fold coproduct of copies of an algebra $A,$ the number of
copies of $A$ in the conclusion, $2^{\cd},$ is, in general,
as large as it can be, embeddability of larger coproducts
being precluded by the size (i.e., cardinality or $\!K\!$-dimension)
of the object we are trying to embed in.
But the assumption that $A$ itself has cardinality or
$\!K\!$-dimension $\leq\cd$ is not forced in that way; we assumed it
so that we could be sure that the coproduct of countably many copies
of $A$ would be embeddable in the object in question.

Can we prove results of the same sort for any larger algebras $A$?

We shall sketch in the next few paragraphs a proof that the symmetric
group $\Sym$ contains a coproduct of two copies of {\em itself}.
Hence, by iteration, it contains coproducts of
all finite numbers of copies of itself, hence, by Lemma~\ref{L.get*w},
a coproduct of countably many copies of itself, hence, by
Theorem~\ref{T.funct}, a coproduct of $2^{\cd}$ copies of itself.
This result, like those that we generalized in preceding sections,
was proved by de~Bruijn in~\cite{debruijn}.
We will then see how to adapt our argument to the case of the monoid
$\Se,$ and, with more work, the associative algebra $\En.$

(In an earlier version of this note, I asked
whether the corresponding result held for the lattice $\Eq.$
An affirmative answer has been given by F.~Wehrung \cite{FW_cP}.)

As indicated above, the hard step, for each of these objects, is
to show that it contains the coproduct of two copies of itself.
Note that to do this for the group $\Sym$
is equivalent to finding two faithful actions of
$\Sym$ on $\Omega$ (or on some set of the same cardinality)
such that there is no nontrivial ``interaction'' between the
permutations giving these actions.
In its most naive form, the idea behind the construction we shall
describe is to take the natural representation of $\Sym$ on $\Omega,$
and the same representation conjugated by a ``random''
permutation $t$ of $\Omega,$ and hope that elements of the two
representations will not interact.

As stated, this is much too naive: no matter how we choose $t$
to eliminate interaction among certain permutations in our two
representations, it will inevitably lead to interaction among others.
However, suppose we replace the set $\Omega$ by the disjoint
union of $\cd$ copies of itself, on each of which
we start with the natural representation of $\Sym,$ and on
each of which we perturb this representation by a different ``$\!t\!$''.
Then we can hope that any given interaction among finitely many
elements of our original and perturbed images of $\Sym$ will be avoided
in at least {\em one} of these copies.
If this is so, then the representation of $\Sym\cP\Sym$ on our union
of copies of $\Omega$ will be faithful.

In particular, we might index our set of copies of $\Omega$ by
the group $\Sym[_{\mathrm{fin}}(\Omega)]$ of all permutations
of $\Omega$ that move only finitely many elements, and on the copy
indexed by each $t$ in that group, let that
$t$ be our perturbing permutation.

There is still one difficulty:  When we construct
a $t$ to prevent interaction in some long expression $w$ in
elements from our two groups, the behavior of $t$ that we need
at one step may be different from the behavior we want at a later step.
To get around this, each copy of $\Omega$ in the above sketch will be
replaced by a disjoint union of countably many copies of itself,
$\Omega\times\omega,$ and $t$ will range over
$\Sym[_{\mathrm{fin}}(\Omega\times\omega)].$
Given a group relation $w=v$ that we want to cause to fail, we will
find that we can select our $t$ and an element
$(p_1,0)\in\Omega\times\{0\}$ so that
as we apply $w$ or $v$ to $(p_1,0),$ that element
is moved by successive occurrences of $t$ into $\Omega\times\{1\},$
$\Omega\times\{2\},$ etc., and on each of
those copies, we shall be able to independently control what $t$ does.

As mentioned, the above technique can also be adapted to the
monoid $\Se,$ and to the $\!K\!$-algebra $\En.$
In the next lemma, the group, monoid, and $\!K\!$-algebra cases are all
stated, and the proof is given for the first two.
I have relegated the longer proof for $\En$ to an appendix,
\S\ref{S.jiggleEn}, so as not to interrupt the flow of the paper.
(Another appendix, \S\ref{S.paths}, gives an alternative construction in
the $\Sym$ case, which I found before encountering de~Bruijn's
papers,
but was not able to adapt to the monoid or $\!K\!$-algebra cases.
It may, however, be of independent group-theoretic interest.)

Recall that we are writing functions to the left of their arguments
(in contrast to the usage in many papers in the theory of infinite
symmetric groups).

\begin{lemma}[{\rm cf.~\cite{debruijn}}]\label{L.double}
\textup{(i)}~ $\Sym$ contains a coproduct of two copies
of itself as a group.\\
\textup{(ii)}~ $\Se$ contains a coproduct of two copies
of itself as a monoid.\\
\textup{(iii)}~ $\En$ contains a coproduct of two copies of itself as
an associative $\!K\!$-algebra.
\end{lemma}\begin{proof}[Proof of \textup{(i)} and \textup{(ii)}]
We shall verify~(ii), then deduce~(i) from it.

Recall that the normal form for an element of the coproduct $M\cP N$
of two monoids is
\begin{xlist}\item\label{x.nmfmg}
$\dots\ \alpha(g_i)\ \beta(g_{i-1})\ \alpha(g_{i-2})\ %
\beta(g_{i-3})\ \dots\ ,$
\end{xlist}
where $\alpha: M\rightarrow M\cP N,$ $\beta: N\rightarrow M\cP N$
are the coprojection maps, the elements $g_k$ with $k$ of one parity
(in~(\ref{x.nmfmg}), the parity of $i)$ are elements of $M-\{1\},$
and those with $k$ of the other parity are elements of $N-\{1\}.$
In~(\ref{x.nmfmg}), I do not explicitly show the first and last factors,
because each may involve either $\alpha$ or $\beta.$
The identity element is given by the empty product~(\ref{x.nmfmg}).

Thus, elements of $\Se\cP\Se$ can be written uniquely
as products~(\ref{x.nmfmg}) in which all $g_i$ come from $\Se-\{1\}.$

To prove~(ii), it suffices to construct a faithful action
of $\Se\cP\Se$ on a set of the same cardinality as $\Omega\<.$
As suggested in the above discussion, that set will be
the disjoint union of a family of copies of $\Omega\times\omega$
indexed by the group $\Sym[_{\mathrm{fin}}(\Omega\times\omega)].$
On every copy of $\Omega\times\omega,$ we let elements
$\alpha(g)$ $(g\in\Se)$ act in the ``natural'' manner,
$g((p,k))=(g(p),k)$ $(p\in\Omega,\ k\in\omega),$ while on the
copy of $\Omega\times\omega$ indexed by
$t\in \Sym[_{\mathrm{fin}}(\Omega\times\omega)],$ we let
$\beta(g)$ act by $t\,g\,t^{-1},$ i.e., the conjugate
by $t$ of that same natural action.

To prove that the resulting action of $\Se\cP\Se$ on our union of copies
of $\Omega\times\omega$ is faithful, assume we are given two distinct
elements of that monoid, say~(\ref{x.nmfmg}) and
\begin{xlist}\item\label{x.nmfmh}
$\dots\ \alpha(h_j)\ \beta(h_{j-1})\ \alpha(h_{j-2})\ \beta(h_{j-3})\ %
\dots\ .$
\end{xlist}
We shall show below how to obtain a
$t\in \Sym[_{\mathrm{fin}}(\Omega\times\omega)]$ such that the induced
actions
of~(\ref{x.nmfmg}) and~(\ref{x.nmfmh}) on $\Omega\times\omega,$ namely
\begin{xlist}\item\label{x.2wt+-}
$\!\dots\ g_i\,\ (t\,g_{i-1}\,t^{-1})\ g_{i-2}\,\ (t\,g_{i-3}\,t^{-1})
\ \dots$\quad and\\[2pt]
$\dots\ h_j\ (t\,h_{j-1}t^{-1})\ h_{j-2}\ (t\,h_{j-3}t^{-1})\ \dots,$
\end{xlist}
act differently on a certain element of $\Omega\times\omega.$

The $t$ we shall construct will be of order~2,
so the above two expressions take the forms
\begin{xlist}\item\label{x.2wt}
$\!\dots\ g_i\,\ t\ g_{i-1}\,\ t\ g_{i-2}\,\ t\ g_{i-3}\,\ t\ \dots$
\quad and\\[2pt]
$\dots\ h_j\ t\ h_{j-1}\ t\ h_{j-2}\ t\ h_{j-3}\ t\ \dots\,.$
\end{xlist}
We may assume, by interchanging~(\ref{x.nmfmg}) and~(\ref{x.nmfmh})
if necessary, that the former expression involves at least as
many factors from $\Sym$ as the latter, and, moreover, that {\em if}
they have the same number of such factors, and have $\!\alpha\!$s
and $\!\beta\!$s in the same places, then for the least value $k$
such that $g_k\neq h_k,$ some element of $\Omega$ on which
$g_k$ and $h_k$ disagree is moved by the former.

If the rightmost term of our original expression~(\ref{x.nmfmg})
is an $\alpha$ term, rather than a $\beta$ term, let us
multiply both lines of~(\ref{x.2wt}) on the right by $t,$
and likewise if the left-hand term of~(\ref{x.nmfmg}) is an $\alpha$
term, let us multiply both lines on the left by $t.$
Since $t$ is going to be invertible,
the non-equality of the new expressions,
which we will prove, is equivalent to the non-equality of the old ones.
The first of the new expressions can now be written more precisely;
the two products have become
\begin{xlist}\item\label{x.2woAt}
$t\ g_n\ t\ \dots\ t\<\ g_i\<\ t\<\ g_{i-1}\<\ t\<\ g_{i-2}\<\ t\<\ %
g_{i-3}\<\ t\ \dots\ t\ g_1\ t,$\quad and\\[2pt]
${}\quad\quad\quad\dots\ t\ h_i\ t\ h_{i-1}\ t\ h_{i-2}\ t\ %
h_{i-3}\ t\ \dots\,.$
\end{xlist}
Clearly, by the assumptions we have made, the first line
of~(\ref{x.2woAt}) has
at least as many occurrences of $t$ as the second.

To construct our promised $t,$
let us now choose, for each $k\in\{1,\dots,n\},$
an element $p_k\in\Omega$ that is moved by $g_k;$ moreover,
we take it to be an element at which
$g_k$ and $h_k$ disagree whenever this is possible; i.e., we require
this for every value of $k$ such that there exists
an element moved by $g_k$ at which $g_k$ and $h_k$ differ.
We then define $t$ to fix all elements of
$\Omega\times\omega$ except the following $2(n+1)$
elements, which we let it transpose in pairs, as shown:
\begin{xlist}\item\label{x.transpose}
$(p_1,0)\leftrightarrow (p_1,1),$\quad
$(g_k(p_k),k)\leftrightarrow (p_{k+1},k{+}1)$ $(1\leq k<n),$
\quad $(g_n(p_n),n)\leftrightarrow (g_n(p_n),n{+}1).$
\end{xlist}
Note that~(\ref{x.transpose}) is consistent:  For $1\leq k\leq n,$
$g_k$ moves $p_k,$ hence the two
elements of $\Omega\times\{k\}$ on which~(\ref{x.transpose})
prescribes (in different ways) the behavior of $t,$ namely $(p_k,k)$
and $(g_k(p_k),k),$ are distinct.

We now see that when the element shown on the first line
of~(\ref{x.2woAt}) is applied to $(p_1,0),$ the successive factors of
that element (reading from the right), namely
$t,$ $g_1,$ $t,$ $g_2,$ $t,\dots\,,$ move it as follows
\begin{xlist}\item\label{x.->->}
$(p_1,0)\mapsto (p_1,1)\mapsto(g_1(p_1),1)\mapsto
(p_2,2)\mapsto(g_2(p_2),2)\mapsto(p_3,3)\mapsto$\\[2pt]
$\quad\quad\dots\mapsto(g_{n-1}(p_{n-1}),n-1)\mapsto (p_n,n)\mapsto
(g_n(p_n),n)\mapsto (g_n(p_n),n+1).$
\end{xlist}
In particular, $(p_1,0)$ is carried from $\Omega\times\{0\}$
into $\Omega\times\{n+1\}.$

When we instead apply the {\em second} line of~(\ref{x.2woAt})
to $(p_1,0)$ there are several possible cases.
If there are fewer factors $h_j$ than $g_i,$ there will be fewer
factors $t$ in that second line than in the first line, so there is no
way the permutation represented by the second line can move an element
from $\Omega\times\{0\}$ into $\Omega\times\{n+1\}.$
If there are the same number, $n,$ of $\!h\!$s as of
$\!g\!$s, but if the $\!\alpha\!$s and $\!\beta\!$s don't appear on
the same factors, then since the first line of~(\ref{x.2woAt})
was adjusted to have a $t$ at each end, the second line will not;
so again there will be fewer factors $t,$ and $(p_1,0)$ cannot
be moved all the way into $\Omega\times\{n+1\}.$

Finally, if there are the same number of factors and
the $\!\alpha\!$s and $\!\beta\!$s appear in the same positions, then
by assumption, for the least $k$ such that $g_k\neq h_k,$
the element $g_k$ moves some element of $\Omega$ at which
these elements disagree (see sentence after~(\ref{x.2wt})), and by
our choice of $p_k,$ the latter will be such an element (first
sentence of paragraph containing~(\ref{x.transpose})).
When we apply the second line of~(\ref{x.2woAt}) to $(p_1,0),$
the input to the factor $h_k$ will be $(p_k,k)$
(since the terms have agreed up to this point),
so the output will be $(h_k(p_k),k)\neq (g_k(p_k),k).$
Thus, our element will fail to be in the unique position
(cf.~(\ref{x.transpose})) from which it can ``catch the boat'' to be
shifted by $t$ from $\Omega\times\{k\}$ to $\Omega\times\{k+1\};$
and since $t$ moves elements by only one level at a time,
our element will not be able to catch up later on.
So the second line of~(\ref{x.2woAt}) does not move $(p_1,0)$
into $\Omega\times\{n+1\},$ hence the two lines
represent distinct elements of $\Se[(\Omega\times\omega)],$
completing the proof of~(ii).

To deduce~(i) from~(ii), note that the normal forms of coproducts of
groups and of monoids are formally the same, hence the inclusion
of $\Sym$ in $\Se$ induces an embedding of $\Sym\cP\Sym$ into
$\Se\cP\Se.$
Since monoid homomorphisms carry invertible elements to
invertible elements, the image of this
copy of $\Sym\cP\Sym$ under the embedding
of statement~(ii) lies in the group $\Sym$ of invertible elements of
$\Se,$ so we have indeed embedded $\Sym\cP\Sym$ in $\Sym,$ as required.

As mentioned earlier, the proof of~(iii) will be given in an
appendix,~\S\ref{S.jiggleEn}.
\end{proof}

By the reasoning sketched at the beginning of this section, we deduce

\begin{theorem}\label{T.bigcP}
\textup{(i)}~ $\Sym$ contains a coproduct of $2^{\cd}$ copies
of itself as a group.\\
\textup{(ii)}~ $\Se$ contains a coproduct of $2^{\cd}$ copies
of itself as a monoid.\\
\textup{(iii)}~ $\En$ contains a coproduct of $2^{\cd}$ copies
of itself as an associative $\!K\!$-algebra.\qed
\end{theorem}

Hence, for instance, if $A$ is any group, not necessarily
of cardinality $\leq\cd,$ that is embeddable
in $\Sym,$ then the coproduct of $2^{\cd}$ copies of $A$
in the category of groups is also embeddable in $\Sym.$

I should mention that at the beginning of this section,
when I said that the cardinal $2^{\cd}$ appearing in
our results was, in general, the best we could hope for, the
phrase ``in general'' was a hedge.
There is an exception, concerning the algebras
$\En$ when $K$ is a field of cardinality $>2^{\cd}.$
We will see at the end of~\S\ref{S.jiggleEn} that in that case,
we can get a stronger conclusion than Theorem~\ref{T.bigcP}(iii).
\vspace{6pt}

Note that Theorem~\ref{T.bigcP}, unlike the results of previous
sections,
says nothing about coproducts in subvarieties of our varieties.
So we ask,

\begin{question}\label{Q.biggerF}
Suppose $A$ is a group, monoid or associative $\!K\!$-algebra
which belongs to a subvariety $\V$ of the variety of
all such algebras, and which is embeddable
in $\Sym,$ $\Se$ or $\En$ respectively.

Must the same be true of the coproduct in $\V$ of two copies of $A$?
\textup{(}If this is indeed true for {\em all} such $A,$ the
corresponding
statement will hold for coproducts in $\V$ of $2^{\cd}$ copies of such
$A,$ by Lemma~\ref{L.get*w} and Theorem~\ref{T.funct}.\textup{)}
\end{question}

In the case where $\V$ is the variety of {\em abelian groups,} or any
subvariety thereof, one has an affirmative answer, for de~Bruijn
\cite[Theorem~4.3]{debruijn_ad} shows that {\em every} abelian group
of cardinality $\leq 2^{\cd}$ is embeddable in $\Sym.$
However, the analog of this stronger statement fails for all varieties
of groups not contained in the variety of abelian groups, by a result
of McKenzie that will be recalled in the next section.

A question similar to the preceding, but concerning
additional constants rather than additional identities, is

\begin{question}\label{Q.amalg}
Suppose $B$ is a subgroup of $\Sym,$ a submonoid of $\Se,$
or a sub-$\!K\!$-algebra of $\En.$

Must $\Sym,$ $\Se$ or $\En$ respectively have a subalgebra containing
$B,$ and isomorphic over $B$ to the coproduct of two copies of
$\Sym,$ $\Se$ or $\En$ with amalgamation of $B$
\textup{(}i.e., isomorphic over $B$ to the pushout, in the variety of
all groups, semigroups, or $\!K\!$-algebras, of the diagram formed
by $B$ and two copies of the indicated algebra; equivalently, to the
coproduct of two copies of that algebra
in the variety of groups, monoids or $\!K\!$-algebras with
distinguished constants corresponding to the elements of~$B)$?

If this is not true in general, does it become true when $B$
has some ``good'' form; e.g., in the case of $\Se,$ when $B$
is a group of invertible elements, or in the case of $\En$ when
$B$ is a division algebra?
\end{question}

Turning back to the argument we used to get
Theorem~\ref{T.bigcP} from Lemma~\ref{L.double},
we should note that a certain fact
was implicitly called on which is true of the varieties of all groups,
all monoids, and all associative $\!K\!$-algebras,
and in many other familiar varieties of algebras, but not in
all -- namely that, given inclusions of algebras $A'\subseteq A$ and
$B'\subseteq B,$ the induced homomorphism of coproducts in our variety,
\begin{xlist}\item\label{x.cPinj}
$A'\cP B'\rightarrow A\cP B,$
\end{xlist}
is also injective.
It is this that allows us to say that if an algebra $A$ contains a
coproduct of two copies of itself, it contains a coproduct of any finite
number of such copies.

An example of a variety $\V$ where the injectivity of
maps~(\ref{x.cPinj})
fails is the variety of groups generated by the infinite dihedral group.
To see this, note that $\V$ satisfies the identity
\begin{xlist}\item\label{x.dih}
$(x^2,y^2)\ =\ 1,$
\end{xlist}
but no identity $x^n=1$ $(n>0).$
Let $A$ and $B$ be infinite cyclic groups $\langle x\rangle$
and $\langle y\rangle;$ these are each free on one generator in $\V.$
Let $A',\ B'$ be the subgroups $\langle x^2\rangle\subseteq A$ and
$\langle y^2\rangle\subseteq B,$ which are isomorphic to $A$ and $B.$
The coproduct $A\cP B$ in $\V$ is the free algebra on $\{x, y\}$
in that variety,
hence is noncommutative, and so the same is true of $A'\cP B'.$
But the image of $A'\cP B'$ in $A\cP B$ is generated by $x^2$ and
$y^2,$ which commute by~(\ref{x.dih}), so the
map $A'\cP B'\rightarrow A\cP B$ is not an embedding.

Though this shows that the principle we used in the proof of
Theorem~\ref{T.bigcP} is not valid in all varieties, it does not show
that the consequence of that principle that we used, concerning objects
containing coproducts of copies of themselves, can fail.
So we ask

\begin{question}\label{Q.2=>3}
Does there exist an algebra $A$ in a variety $\V$ such
that $A$ contains the coproduct in $\V$ of two copies of itself,
but not the coproduct of three such copies?
\end{question}

The preceding results about fitting into $\Sym,$ $\Se$ and $\En$
multiple copies of themselves suggest questions about fitting
these objects into each other, in various ways.
We record two easy results in this direction:

\begin{lemma}\label{L.eachother}
\textup{(i)}~ $\En$ contains an embedded copy of
$K\,\Se,$ the monoid algebra over $K$ on the monoid $\Se.$\\
\textup{(ii)}~ $\Se$ contains an embedded copy of $\Eq_\wedge,$
i.e., $\Eq$ made a monoid under the meet operation~$\wedge.$
\end{lemma}\begin{proof}
(i)~ Let $K\langle\Omega\rangle$ denote the free associative
$\!K\!$-algebra on $\Omega\<.$
Then the action of $\Se$ on $\Omega$ induces an action of $\Se$
on $K\langle\Omega\rangle$ by $\!K\!$-algebra endomorphisms,
which extends to an action of the monoid algebra  $K\,\Se$
by vector-space endomorphisms of $K\langle\Omega\rangle.$
Since $K\langle\Omega\rangle$ has the same vector-space
dimension, $\cd,$ as $V,$ our assertion will follow if we can show
that the endomorphisms of $K\langle\Omega\rangle$ induced by any finite
family $g_1,\dots,g_n$ of distinct elements of $\Se$
are $\!K\!$-linearly independent.

Given such $g_1,\dots,g_n,$ let us choose $p_1,\dots,p_m\in\Omega$ such
that no two of $g_1,\dots,g_n$ behave the same on all of these elements.
Regarding $p_1,\dots,p_m$ as members of the free generating
set $\Omega$ of $K\langle\Omega\rangle,$ we can form the product
$p_1\dots p_m$ therein, and observe that the actions of
$g_1,\dots,g_n$ take this monomial to distinct monomials,
hence are indeed $\!K\!$-linearly independent.

(ii)~ If we regard equivalence relations on $\Omega$
as subsets of $\Omega\times\Omega,$ then the meet operation on $\Eq$
is the restriction of the
intersection operation on $\Pw[\Omega\times\Omega],$ hence it
will suffice to embed $\Pw[\Omega\times\Omega]_\cap$ in $\Se.$
Since $\cd[\Omega\times\Omega]=\cd=\cd[\Omega\times 2],$ we can
do this, in turn, if we can embed $\Pw_\cap$ in $\Se[(\Omega\times 2)].$
To do this, let us send each $S\subseteq\Omega$
to the endomap of $\Omega\times 2$ that fixes all elements $(p,0),$
and also all elements $(p,1)$ with $p\in S,$ but sends
$(p,1)$ to $(p,0)$ if $p\notin S.$
The verification that this is a monoid homomorphism, and indeed
an embedding, is straightforward.
\end{proof}

We will see in \S\ref{S.notEq} that the analog of
statement~(ii) above with ``meet'' replaced by ``join'' is false.

It is also interesting to note that the analog of~(i) fails if
$\Omega$ replaced by a finite set with $n\geq 2$ elements.
Indeed, for $n\geq 3,$
even the group algebra $K\,\Sym$ cannot be embedded in $\En,$
for it has dimension $n!,$ while $\En$ has dimension only $n^2.$
In particular, the $n!$ permutation matrices do not generate a copy
of the group algebra -- they are not linearly independent.
To get the nonembeddability
statement for $K\,\Se$ when $n=2,$ let $R=K\,\Se[(2)],$ and
note by comparing dimensions that a $\!K\!$-algebra embedding of $R$
in $M_2(K)$ would have to be an isomorphism.
Let $z\in\Se[(2)]$ be the map taking both elements of
$2$ to $0,$ and note that it satisfies the left-zero
identity $(\forall\,a)\ za=z.$
Hence $zR\subseteq Kz,$ hence $z\,R\,(1-z)=0;$ but $M_2(K)$ has no
idempotent $z$ with this property other than $0$ and $1.$

\section{Restrictions on groups embeddable in $\Sym$ and
monoids embeddable in $\Se.$}\label{S.notSym}

With such vast classes of groups, monoids, associative
$\!K\!$-algebras and lattices embeddable in $\Sym,$ $\Se,$ $\En,$ and
$\Eq,$ it is natural to ask whether there are
groups, etc., of cardinality, respectively $\!K\!$-dimension,
$\leq 2^{\cd},$ that are not so embeddable.

For the case of groups, de~Bruijn~\cite{debruijn} showed, in effect,
that for any set $I$ of cardinality $>\cd,$
the group presented by generators $x_i$ $(i\in I)$ and relations
\begin{xlist}\item\label{x.235:2}
$x_i^2=1$ $(i\in I),$
\end{xlist}
\begin{xlist}\item\label{x.235:3}
$(x_i\<x_j)^3=1$ $(i,j\in I,$ distinct),
\end{xlist}
\begin{xlist}\item\label{x.235:5}
$(x_i\<x_j\<x_k\<x_l)^5=1$ $(i,j,k,l\in I,$ distinct),
\end{xlist}
cannot be embedded in $\Sym.$
Note that the fact that the indices in~(\ref{x.235:3})
and~(\ref{x.235:5}) are required to be distinct keeps this
system of groups from having the form to which the results
of~\S\ref{S.cP+F} apply.
(If those indices were not required to be distinct,
then setting $k=i,\ l=j$ in~(\ref{x.235:5}) would give
$(x_i x_j)^{10}=1,$ which, combined with~(\ref{x.235:3}),
would give $x_i x_j=1,$ making the group collapse to $Z_2.)$

On the other hand, de~Bruijn claimed in~\cite{debruijn}
that his result corresponding to Theorem~\ref{T.funct} showed that the
restricted direct product (called in~\cite{debruijn} the direct product)
of $2^{\cd}$ copies of any group $A$ of cardinality
$\leq\cd$ could be embedded in $\Sym$ -- not noticing that because
the commutativity
relations which the restricted direct product construction imposes
on elements of different copies of $A$ fail to hold among
elements of a single copy (unless $A$ is commutative), that
result is not applicable.
In~\cite{debruijn_ad} he corrected this error,
noting that the argument is only valid when $A$ is abelian,
and posed his earlier assertion as an open question.

That question was answered in the negative by McKenzie~\cite{mckenzie},
who showed that if $G$ is a group such that for some index set $I$ with
$\cd[I]>\cd$ there are elements $x_i, y_i\in G$ $(i\in I)$ satisfying
\begin{xlist}\item\label{x.cmxcm}
$x_iy_j=y_jx_i$
whenever $i\neq j,$ but $x_iy_i\neq y_ix_i,$
\end{xlist}
then $G$ cannot be embedded in $\Sym.$

We shall see that McKenzie's criterion is
an instance of more general facts.
By a {\em centralizer subgroup} in a group $G,$ let us understand a
subgroup of the form
\begin{xlist}\item\label{x.C_G}
$C_G(X)\ =\ \{g\in G\mid(\forall\,x\in X)\ gx=xg\}$
\end{xlist}
for some subset $X\subseteq G.$
A subgroup $H<G$ is clearly a centralizer subgroup
if and only if $H=C_G(C_G(H)).$
Recall also that a {\em jump} in a totally ordered set means
a pair of elements $x<y$ such that $\{z\mid x<z<y\}$ is empty.
A totally ordered set without jumps can have subsets with jumps;
for instance, the set of reals or of rationals has none, but
their subset $\Z$ has countably many.
We shall see below that the lattice of centralizer subgroups of
$\Sym,$ and hence of any group embeddable therein,
contains no chains with $>\cd$ jumps, while
the lattice of centralizer subgroups of a group with a family of
elements satisfying~(\ref{x.cmxcm}) does have such chains.

In fact, we shall prove the former result not only for
centralizer subgroups, but for subsets of $\Sym$ defined by
arbitrary systems of equations (in several variables)
with constants in $\Sym,$ which will also
yield a quick proof of de~Bruijn's example.
Our result will follow from the fact that such solution subsets are
closed in the function topology on $\Sym,$
together with the following lemma in general topology.

Note that unless explicitly stated,
we do not assume topologies to be Hausdorff.
The function topology, to which we will apply the lemma in this
section, is Hausdorff, but in the next section we will apply the
same lemma to both Hausdorff and non-Hausdorff topologies.

\begin{lemma}\label{L.cdB}
Let $T$ be a topological space having an infinite basis \textup{(}or
more generally, subbasis\textup{)} $B$ of open sets.
Then the lattice of open subsets of $T$ contains no chain
with $>\cd[B]$ jumps.
Hence its opposite, the lattice
of closed subsets of $T,$ also has no such chains.
In particular, that lattice contains no well-ordered
or reverse-well-ordered chains of cardinality $>\cd[B].$
\end{lemma}\begin{proof}
The case where $B$ is a subbasis reduces to that in
which it is a basis, since in the former case, a basis
is given by the set of intersections of finite subsets of $B,$
and for $B$ infinite there are only $\cd[B]$ of these.
So we assume $B$ a basis.

Suppose $C$ is a chain of open subsets of $T.$
For each jump $U\subset V$ in $C,$ let us choose a point
$p_{(U,V)}\in V-U.$
Since $V$ is a neighborhood of $p_{(U,V)},$ our basis
$B$ contains some
subneighborhood $N_{(U,V)}\subseteq V$ of $p_{(U,V)};$
let $N_{(U,V)}$ be so chosen for each jump $U\subset V.$
Then if $U\subset V$ and $U'\subset V'$ are distinct jumps, say with
$U\subset V\subseteq U'\subset V',$ we must have
$N_{(U,V)}\neq N_{(U',V')},$ since $N_{(U,V)}\subseteq V,$ while
$N_{(U',V')}$ contains $p_{(U',V')}\notin U'.$
Hence distinct jumps in $C$ give distinct elements $N_{(U,V)}\in B,$
so the number of jumps in $C$ does not exceed $\cd[B].$
\end{proof}

Recall next that if $\Omega$ and $\Omega'$ are sets, and we give
$\Omega'$
the discrete topology, then the {\em function topology} on the set
of all maps $\Omega\rightarrow \Omega'$ has for a subbasis of open sets
the sets $U_{y,x}=\{f:\Omega\rightarrow \Omega'\mid f(x)=y\}$
$(x\in\Omega,\ y\in \Omega'),$ since a basis for the open sets
of $\Omega'$ is given by the singletons $\{y\}.$
In particular, if $\Omega$ is infinite, the function topology on the
set $\Se$ has a subbasis of cardinality $\cd[\Omega\times\Omega]=\cd.$
Hence for any set $J,$ the direct product of a $\!J\!$-tuple of copies
of this space has a subbasis of cardinality $\cd\,\cd[J].$
So we get

\begin{corollary}\label{C.funcbss}
If $\cd[J]\leq\cd,$ then the lattice of subsets of $\Se^J$ closed in
the product topology on that set induced by the function topologies
on the factors $\Se$ has no chains with
\linebreak[3] $>\nolinebreak\cd$ jumps.\qed
\end{corollary}

Let us now connect this topology with our algebraic structure.
It is straightforward to verify that the operation of composition
on $\Se$ is continuous in the function topology;
moreover, on its subset $\Sym,$ the operation
of functional inverse is also continuous, since
it simply interchanges $U_{y,x}\cap\Sym$ and $U_{x,y}\cap\Sym$
for all $x$ and $y.$
(We remark, however, that $\Sym$ is not closed in $\Se.)$
Hence given a pair of monoid words (respectively
group words) $v,\ w,$ in a variable $t$
and constants from $\Se$ (respectively, from $\Sym),$
the solution set $\{a\mid v(a)=w(a)\}$ will be closed in the
function topology on $\Se$ (respectively, $\Sym).$
More generally, we may look at the solution set
of any family of such pairs of words in any number of variables.
Let us set up notation for such sets in an arbitrary algebra.
\begin{definition}\label{D.S^=}
For $A$ an algebra in a variety $\V,$ and $J$ a set, we shall
understand a {\em principal solution set} in $A^J$ to
mean a set of the form
\begin{xlist}\item\label{x.solset}
$S_{v=w}\ =\ \{a=(a_j)_{j\in J}\in A^J\mid v(a)=w(a)\}\ \subseteq\ A^J,$
\end{xlist}
where $v$ and $w$ are words in a $\!J\!$-tuple of variables
$(t_j)_{j\in J},$ constants from $A,$ and the operations of $\V.$

A {\em solution set} in $A^J$ will mean the intersection of an arbitrary
family of principal solution sets.
We shall denote by $L^{\strt=}_{A,\,J}$ the complete lattice of
all solution sets in $A^J.$
\end{definition}

Here we understand the intersection of the {\em empty} family of
principal solution sets to be the whole set $A^J.$
Thus, $L^{\strt=}_{A,\,J}$ is, as asserted, a complete lattice,
the join of any $X\subseteq L^{\strt=}_{A,\,J}$ being the
intersection of those
principal solution sets that contain all members of $X.$
(The superscript ``$\!=\!$'' in $L^{\strt=}_{A,\,J}$ indicates that our
solution sets are defined by equations, as in~(\ref{x.solset}).
In the next section we shall also make use
of solution sets defined by inequalities.)

In the next result, though $J$ is allowed to have cardinality up to
$\cd,$ the most common cardinality in our applications will be $1.$

\begin{theorem}\label{T.gp+md.crit}
Let $J$ be any set of cardinality $\leq\cd.$
Then $L^{\strt=}_{\Sym,\,J}$ contains no chains with
$>\nolinebreak\cd$ jumps.
Hence the same is true of $L^{\strt=}_{G,\,J}$
for any group $G$ embeddable in $\Sym.$

Likewise, $L^{\strt=}_{\Se,\,J}$ contains no chains with $>\cd$ jumps;
hence the same is true of $L^{\strt=}_{M,\,J}$
for any monoid $M$ embeddable in $\Se.$

In particular, for $G$ a group embeddable in $\Sym$ or $M$ a monoid
embeddable in $\Se,$ the lattice $L^{\strt=}_{G,\,J},$
respectively $L^{\strt=}_{M,\,J}$ contains no well-ordered or
reverse-well-ordered chain of cardinality $>\cd.$
\end{theorem} \begin{proof}
The assertions about chains
in $L^{\strt=}_{\Sym,\,J}$ and $L^{\strt=}_{\Se,\,J}$ are clear
from Corollary~\ref{C.funcbss} and the continuity of our operations;
it remains to deduce the corresponding statements for objects embeddable
in $\Sym$ and $\Se.$

If $G$ is a group embeddable in $\Sym,$ let us assume for notational
convenience that it is a subgroup, and map $L^{\strt=}_{G,\,J}$
to $L^{\strt=}_{\Sym,\,J}$ by sending each solution set $S$ in $G^J$
to the solution set in $\Sym^J$ of the set of all equations
(in a $\!J\!$-tuple of variables, with constants in $G)$ that
are satisfied on $S.$
This map is easily seen to be an embedding of partially ordered
sets, hence the result on chains in $L^{\strt=}_{\Sym,\,J}$ implies
the same conclusion for chains in $L^{\strt=}_{G,\,J}.$
The same argument works for monoids.
\end{proof}

Some observations on the above proof:
Given groups $G<H,$ one cannot embed $L^{\strt=}_{G,\,J}$ in
$L^{\strt=}_{H,\,J}$ by simply sending the solution set of every system
of equations in $G$ to the solution set of the same system in $H.$
This does not give a well-defined function,
since equations over $G$ having the same solution set
in $G^J$ may have different solution sets in $H^J.$
(Consider, for instance, centralizer subgroups of various sets
in an abelian group, and of the same sets in a nonabelian overgroup.)
The construction of the above proof does give an
order-embedding of $L^{\strt=}_{G,\,J}$ into $L^{\strt=}_{H,\,J},$
but in general this respects neither meets nor joins;
the former because the set $H^J$ is larger than $G^J;$ the latter
because the set of equations with constants in $H$ is larger
than the set of equations with constants in $G.$
Another order-embedding of $L^{\strt=}_{G,\,J}$ in $L^{\strt=}_{H,\,J}$
is gotten by sending every $S\in L^{\strt=}_{G,\,J}$
to the set of elements of $H^J$ satisfying all equations
{\em with constants in $H$} satisfied on $S;$
it also respects neither meets nor joins, in general.
\vspace{6pt}

For our first application of Theorem~\ref{T.gp+md.crit}, note that
the centralizer subgroups~(\ref{x.C_G}) in a group $G$ form a complete
lattice, which as a partially ordered set (and indeed, as a complete
lower semilattice) is embedded in $L^{\strt=}_{G,1}.$
Hence we have

\begin{corollary}\label{C.cntrlzr}
No group having a chain of centralizer
subgroups with $>\cd$ jumps is embeddable in $\Sym.$
In particular \textup{(McKenzie, \cite{mckenzie})},
if a group $G$ contains, for some set $I$ with $\cd[I]>\cd,$
elements $x_i,\,y_i$ $(i\in I)$ satisfying\textup{~(\ref{x.cmxcm})},
then $G$ is not embeddable in $\Sym.$
\end{corollary}\begin{proof}
The first statement is clear from the first paragraph
of Theorem~\ref{T.gp+md.crit}.
In the situation of the second statement, we may, by reindexing,
assume $I$ to be a cardinal $\kappa>\cd.$
For each $\alpha\in\kappa,$ let
$X_\alpha=\{x_\beta\mid \beta>\alpha\}.$
The $X_\alpha$ form a descending chain of subsets,
hence their centralizers $C_G(X_\alpha)$ form
an ascending chain of centralizer subgroups.
Note that each $C_G(X_\alpha)$ contains those elements $y_\gamma$ with
$\gamma\leq\alpha$ and no other $y_\gamma,$
hence the $C_G(X_\alpha)$ are distinct.
Thus we have a well-ordered chain of centralizer subgroups of
cardinality $\kappa>\cd;$ hence $G$ is not embeddable in $\Sym.$
\end{proof}

The cardinality conditions in the above result
are sharp:  If we take any set $I$
of cardinality $\leq\cd,$ and any $\!I\!$-tuple $G_i$ of
nonabelian groups each of cardinality $\leq\cd,$ then
their restricted direct product has cardinality $\leq\cd,$
hence is embeddable in $\Sym,$ though it
contains elements $x_i,\,y_i$ satisfying~(\ref{x.cmxcm}).

Turning to de~Bruijn's relations~(\ref{x.235:2})-(\ref{x.235:5}),
note that for any index set $I$ of cardinality
$\leq\cd,$ if we take an element $p_0\in\Omega$ and
an $\!I\!$-tuple of elements $p_i\in\Omega$ distinct from $p_0$ and
from each other, and for each $i\in I$ let $x_i\in\Sym$
be the transposition that interchanges $p_0$ and $p_i$ and fixes
all other elements, then any product of $n$ distinct
elements  $x_i$ is an $\!(n+1)\!$-cycle in $\Sym,$
from which we see that the $\!I\!$-tuple
$(x_i)_{i\in I}$ satisfies~(\ref{x.235:2})-(\ref{x.235:5}).
This shows that the cardinality conditions in the
next result are sharp.
That result, in fact, does without~(\ref{x.235:2}), at
the small price of adding two sentences at the start of the proof.

\begin{corollary}[{\rm cf.~\cite[Theorem~5.1]{debruijn}}]\label{C.235}
No group $G$ containing a family $(x_i)_{i\in I}$ of distinct elements
satisfying~\textup{(\ref{x.235:3})} and~\textup{(\ref{x.235:5})},
where $\cd[I]>\cd,$ is embeddable in $\Sym.$
\end{corollary}\begin{proof}
If there are any pairs $i\neq i'\in I$
such that $x_i$ and $x_{i'}$ are inverse to one another,
then dropping one member of each such pair does not decrease $\cd[I].$
Hence we may assume there are no such pairs.

Let us also assume, by reindexing, that $I$ has the
form $\kappa\times 2,$ where $\kappa$ is a cardinal $>\cd,$
so that our given elements have the form $x_{\beta,i}$
$(\beta\in\kappa,\ i=0,1).$
For each $\alpha\in\kappa,$ let us define the solution set
$S_\alpha=\{(y,z)\in G^2\mid(\forall\,\beta>\alpha)\ %
(x_{\beta,0}\ x_{\beta,1}\ y\ z)^5=1\}.$
By~(\ref{x.235:5}), this set contains the pair
$(x_{\gamma,0},x_{\gamma,1})$ whenever $\gamma\leq\alpha.$

However, it contains no pair
$(x_{\beta,0},x_{\beta,1})$ with $\beta>\alpha.$
Indeed, if it did, we would have
\begin{xlist}\item\label{x.b0b1b0b1}
$1\,=\,(x_{\beta,0}\ x_{\beta,1}\ \linebreak[0]
x_{\beta,0}\ x_{\beta,1})^5\,=\,(x_{\beta,0}\ x_{\beta,1})^{10}.$
\end{xlist}
But by~(\ref{x.235:3}), $(x_{\beta,0}\ x_{\beta,1})^3=1.$
Combining these equations we get $x_{\beta,0}\ x_{\beta,1}=1,$
contradicting our assumption that for $i$ and $i'$ distinct,
$x_i$ and $x_{i'}$ are not inverses.

Hence, the sets $S_\alpha$ are distinct, and so form a well-ordered
chain of cardinality $\kappa$ in $L^{\strt=}_{G,\,2},$ from which
nonembeddability of $G$ in $\Sym$ follows by Theorem~\ref{T.gp+md.crit}.
\end{proof}

What about applications of Theorem~\ref{T.gp+md.crit} to monoids?
Well, the monoid homomorphisms from a group $G$ to
a monoid $M$ are the group homomorphisms from $G$ to the group
of invertible elements of $M;$ hence the above two corollaries can
also be viewed as giving monoids that are not embeddable in $\Se.$

Here is a more genuinely monoid-theoretic application.
Consider again an element $p_0\in\Omega$
and a family of distinct elements $p_i\in\Omega-\{p_0\}$ indexed by
a set $I$ of the same cardinality as $\Omega\<.$
Let $y\in\Se$ be the map sending all elements to $p_0,$ while
for each $i\in I,$ let $x_i$ be the map sending everything
except $p_i$ to $p_0,$ and fixing $p_i.$
Then we see that
\begin{xlist}\item\label{x.xxy}
for $i,j\in I,$ $x_i x_j = y$ if and only if $i\neq j.$
\end{xlist}
This gives $\cd$ such elements $x_i;$
but an application of Theorem~\ref{T.gp+md.crit}, following the
same pattern as the two preceding results, shows that
we cannot get a family of $>\cd$ such elements; hence

\begin{corollary}\label{C.eqprod}
No monoid containing an element $y,$ and a family of elements $x_i,$
distinct from $y,$ indexed by a set $I$ of cardinality $>\cd,$ and
satisfying~\textup{(\ref{x.xxy})}, is embeddable in $\Se.$\qed
\end{corollary}

In our examples of nonembeddability using the
group conditions (\ref{x.235:3})-(\ref{x.235:5}) and~(\ref{x.cmxcm}),
and the monoid conditions~(\ref{x.xxy}), we could have
asserted much more than the existence of a chain with $\cd[I]$ jumps.
For example, given an $\!I\!$-tuple of elements
satisfying~(\ref{x.cmxcm}), distinct subsets of $\{x_i\mid i\in I\}$
have centralizers containing distinct subsets of the $y_i$ (indexed
by the complementary subsets of $I),$ so we in fact get a copy
of the whole partially ordered set $\Pw[I]$ in $L^{\strt=}_{G,1};$
and the corresponding observations hold for the other two examples.

However, there are examples that give large chains of solution sets
without (as far as I can see) giving so much more as well.
If we take the group or monoid presented by~(\ref{x.xyijR}),
respectively~(\ref{x.xiR}), with $I$ a cardinal $\kappa,$ then it will
have a chain of centralizers, respectively fixed sets,
order-isomorphic to $\kappa$ (as well as one
of the opposite order type), but there is no apparent reason why
it should have, say, any large antichain of solution sets.

It is interesting that while the constructions
of~(\ref{x.xyijR}) and~(\ref{x.xiR}) with $I=\R$ give, as
we saw earlier, groups and monoids embeddable in $\Sym[(\omega)],$
respectively $\Se[(\omega)],$ the above paragraph shows that the
contrary is
true for the same constructions with $I=\cd[\R]$ as a well-ordered set.
We likewise get nonembeddability when $I=\R\times 2$ under
lexicographic order (since it also
has uncountably many jumps), and for $\R\times\R$ under
lexicographic order (since this contains the preceding ordered set).

Here is another application of our observations on~(\ref{x.xyijR}).
For every real number $c,$ let $G_c$ be the group presented
by generators $x_r,\,y_r$ $(r\in\R)$ and relations
\begin{xlist}\item\label{x.xycR}
$x_r\<y_s = y_s\<x_r$ for all $r,\,s\in\R$ such that $s\geq r+c.$
\end{xlist}
Clearly, each $G_c$ is isomorphic to $G_0,$ by an isomorphism that
fixes the $x_r$ and takes each $y_s$ to $y_{s-c}.$
But we have noted that (by Theorem~\ref{T.F(R)}) $G_0$ is
embeddable in $\Sym[(\aleph_0)];$ hence so is every $G_c.$
Now let $G_{0^+}$ denote the group with the
same generators, but having for relations the union of the sets
of relations defining $G_c$ for all $c>0;$ in other words,
\begin{xlist}\item\label{x.>0}
$x_r\<y_s = y_s\<x_r$ for all $r,\,s\in\R$ with $s>r,$
\end{xlist}
and let us define therein, for every $c\in\R,$ the centralizer subgroups
\begin{xlist}\item\label{x.ScSc+}
$S_c\ =\,\ \{g\in G_{0^+}\mid (\forall\,s\geq c)\ g\<y_s=y_s\<g\},$
\\[2pt]
${}\,S_{c^+}=\ \{g\in G_{0^+}\mid (\forall\,s>c)\ g\<y_s=y_s\<g\}.$
\end{xlist}
I claim that these form a chain, with jumps $S_c\subset S_{c^+},$
and with inclusions $S_{c^+}\subset S_d$ whenever $c<d.$
Indeed, $S_c\neq S_{c^+}$ because $x_c\in S_{c^+}-S_c,$ and similarly
$S_{c^+}\neq S_d$ by considering $x_e$ for any $e$ with $c<e<d.$
This chain is isomorphic to $\R\times 2,$ so $G_{0^+}$ is not
embeddable in $\Sym[(\aleph_0)],$ though it is a direct limit,
via surjective homomorphisms,
of the groups $G_c$ $(c>0),$ which are so embeddable.

Though a positive answer seems implausible, let us ask

\begin{question}\label{Q.gps}
Is the criterion of the first paragraph of Theorem~\ref{T.gp+md.crit}
also {\em sufficient} for a group of
cardinality $\leq 2^{\cd}$ to be embeddable in $\Sym$?
If so, is it sufficient that it hold for all finite $J$?
For $J=1$?
\end{question}

In an earlier version of this note, I asked the same questions
for $\Se;$ but the possibility of an affirmative answer is
now precluded by a result of Wehrung~\cite{FW_op},
showing that $\Se\op$ cannot be embedded in $\Se;$
indeed, since the lattices $L^{\strt=}_{\Se,\,J}$ and
$L^{\strt=}_{\Se\op,\,J}$ are isomorphic, that result shows that
no condition on the lattices $L^{\strt=}_{M,\,J}$ can be equivalent
to embeddability of $M$ in $\Se.$
Conceivably, however, one could construct lattices of solution sets
not using all the sets $S_{v=w},$ but some subfamilies of these that are
not invariant under reversing the orders of factors in the words $u$
and $v,$ such that conditions on these solution sets would characterize
embeddability in $\Se.$

It is easy to {\em formally} strengthen
Theorem~\ref{T.gp+md.crit} in several ways.
First, since the sets $S_{v=w}$ of~(\ref{x.solset})
are closed in the function topology, so are finite
unions of such sets, which we might write
\begin{xlist}\item\label{x.usolset}
$S_{(v_1=w_1)\vee\dots\vee(v_n=w_n)}\ =
\ S_{v_1=w_1}\cup\dots\cup S_{v_n=w_n}.$
\end{xlist}
So if we let $L^{=,\vee}_{\Sym,\,J},$
denote the lattice of arbitrary intersections
of families of finite unions~(\ref{x.usolset}), these will also satisfy
the conditions on chains given by Theorem~\ref{T.gp+md.crit}.

Secondly, the conclusion of Theorem~\ref{T.gp+md.crit} only states
one particular consequence of embeddability of our lattice of
solution sets in the lattice of closed sets
of a topology generated by $\leq\cd$ elements.
We will examine the latter condition further in an appendix,
\S\ref{S.=>=>}.
Meanwhile, we ask

\begin{question}\label{Q.solset+}
If the criteria of Theorem~\ref{T.gp+md.crit} are not
sufficient for a group $X$ of cardinality $\leq 2^{\cd}$
to be embeddable in $\Sym,$ do they
become so if we replace the lattices $L^{\strt=}_{X,\,J}$ of that
theorem by the larger lattices $L^{=,\vee}_{X,\,J},$ and/or
strengthen the condition on jumps in chains to the
condition that our lattice of solution sets be
embeddable as a partially ordered
set in the lattice of closed subsets of a topological space
with a basis of cardinality $\leq\cd$
\textup{(}cf.\ \S\ref{S.=>=>}\textup{)}?

Do there, at least, exist groups whose embeddability in $\Sym$ is
precluded by one of these strengthened conditions, but not
by the conditions of Theorem~\ref{T.gp+md.crit}?
Here we may ask the same question for embeddability of monoids in $\Se.$
\end{question}

We mentioned (following Question~\ref{Q.biggerF})
de~Bruijn's result that every abelian group of
cardinality $\leq 2^{\cd}$ is embeddable in $\Sym.$
However, not every commutative monoid of cardinality $\leq 2^{\cd}$
embeds in $\Se:$
the presentations~(\ref{x.xiR}) give commutative monoids,
but we saw that for $I$ a cardinal $>\nolinebreak\cd,$
the resulting monoid is not so embeddable.
So we ask

\begin{question}\label{Q.VinSe}
Which varieties $\V$ of monoids have the property that every monoid
in $\V$ of cardinality $\leq 2^{\cd}$ is embeddable in $\Se$?
\end{question}

By de~Bruijn's result on abelian groups, this is true of every
variety of commutative monoids satisfying an identity $x^n=1,$ since
such monoids are essentially abelian groups of exponent $n,$ hence
embeddable in $\Sym$ as groups.
I don't know any other examples.

Returning to groups, suppose we write
$\Sym[_<(\Omega)]\subset\Sym$ for the normal subgroup of permutations
that move fewer than $\cd$ elements.
De~Bruijn \cite[Theorem~4.4]{debruijn} showed that $\Sym$ could
be embedded in $\Sym/\Sym[_<(\Omega)],$ while
McKenzie~\cite[Corollary~3]{mckenzie} showed that
$\Sym/\Sym[_<(\Omega)]$ contains a restricted direct product of $>\cd$
copies of itself, and hence, by Corollary~\ref{C.cntrlzr}, cannot
be embedded in $\Sym.$

\begin{question}\label{Q.S/fin}
What restrictions does embeddability in
$\Sym/\Sym[_<(\Omega)]$ imply for a group of cardinality $\leq 2^{\cd}$?
\end{question}

Under the assumption of the General Continuum Hypothesis,
J\'onsson~\cite{jonsson} shows, inter alia, that for every %yes!
uncountable
cardinal $\kappa$ there exists a group of cardinality $\kappa$ which
contains isomorphic copies of {\em all} groups of cardinality $\kappa.$
We have seen that $\Sym$ is {\em not} such a group for
$\kappa=2^{\cd}.$
Felgner and Haug~\cite{fe+ha} show that under
certain set-theoretic hypothesis, neither is $\Sym/\Sym[_<(\Omega)].$

\section{Restrictions on lattices embeddable in $\Eq.$}\label{S.notEq}

If we want to adapt the technique of the preceding section to
get restrictions on lattices embeddable in the lattice $\Eq$
of equivalence relations on $\Omega,$ we must decide what
topology on that lattice to use in place of the function topology.
One approach is to regard binary relations on $\Omega$
as elements of $\Pw[\Omega\times\Omega]=2^{\Omega\times\Omega},$ i.e.,
as functions $\Omega\times\Omega\rightarrow 2=\{0,1\},$ and use the
function topology on that set induced by the discrete topology on $2.$
A subbasis of open subsets of $\Pw[\Omega\times\Omega]$ under
this topology is given by the sets
\begin{xlist}\item\label{x.Upqc}
$U_{p,q}=\{R\in\Pw[\Omega\times\Omega]\mid (p,q)\in R\}$~ and\\[2pt]
$^c U_{p,q}=\{R\in\Pw[\Omega\times\Omega]\mid (p,q)\notin R\}.$
\end{xlist}
We see that each of these sets is clopen (closed and open) in the
topology so defined, and that this subbasis has cardinality $\cd.$

By abuse of notation, in speaking of subsets of $\Eq$
let us write $U_{p,q}$ for $U_{p,q}\cap\Eq$
and $^c U_{p,q}$ for $^c U_{p,q}\cap\Eq$
(just as, in introducing the function topology on $\Se,$ we earlier
wrote $U_{p,q}$ for what we would now describe as $U_{p,q}\cap\Se).$

In this topology, one finds that the meet operation,
i.e., intersection as subsets of $\Omega\times\Omega,$ is continuous,
but that the join operation is not.
To see the first fact, note that under the map
\begin{xlist}\item\label{x.meet}
$\wedge:\ \Eq\times\Eq\rightarrow \Eq,$
\end{xlist}
the inverse image of every $U_{p,q}$ is
the open rectangle $U_{p,q}\times U_{p,q},$ while the inverse
image of $^c U_{p,q}$ is the union
$(^c U_{p,q}\times\Eq)\cup(\Eq\times\,^c U_{p,q}),$ and both
these sets are open.
Under the operation
\begin{xlist}\item\label{x.join}
$\vee:\ \Eq\times\Eq\rightarrow \Eq,$
\end{xlist}
the inverse image of $U_{p,q}$ is still open: it
is an infinite union of finite intersections of sets of the forms
$U_{r,s}\times\Eq$ and $\Eq\times U_{r,s},$ one such intersection for
each finite chain of relations which, if they hold in a pair of
equivalence relations on $\Omega,$ witness the conclusion that
$(p,q)$ belongs to the join of those equivalence relations.
But the inverse image of $^c U_{p,q}$ becomes, by the same reasoning,
an infinite {\em intersection} of clopen sets, which will not be open.
(Essentially because no {\em finite} set of relations and
negations of relations can witness the absence of
$(p,q)$ from the join of two equivalence relations.)

There are two ways to respond to this difficulty.
One is to formulate a criterion in terms of words in the meet
operation alone.
This corresponds to regarding $\Eq$ as a meet-semilattice, $\Eq_\wedge.$
In the notation of Definition~\ref{D.S^=},
we then look at the lattices $L^{\strt=}_{\Eq_\wedge,\,J},$ and get

\begin{theorem}\label{T.semilatcrit}
Let $J$ be any set of cardinality $\leq\cd.$
Then $L^{\strt=}_{\Eq_\wedge,\,J}$ contains
no chains with $>\cd$ jumps.
Hence for any lower semilattice $A$ embeddable in $\Eq_\wedge,$ the
lattice $L^{\strt=}_{A,\,J}$ contains no chains with $>\cd$ jumps.
In particular, for any lattice $A$ embeddable in $\Eq,$
the lattice $L^{\strt=}_{A_\wedge,\,J}$ has no such chains.
\end{theorem}\begin{proof}
This can be gotten by the same method as Theorem~\ref{T.gp+md.crit}, or
deduced therefrom using Lemma~\ref{L.eachother}(ii).
\end{proof}

The following application of this result shows that the meet-join
asymmetry we have come up against is real;
and the second statement gives the promised example showing
that the analog of Lemma~\ref{L.eachother}(ii)
with $\Eq_\vee$ in place of $\Eq_\wedge$ is false.

\begin{proposition}\label{P.mtvsjn}
The largest cardinality of a set $I$ such that $\Eq$ contains an
element $z$ and elements $x_i\neq z$ $(i\in I)$ satisfying
\begin{xlist}\item\label{x.meetz}
$x_i\wedge x_j=z$ for all distinct $i,j\in I$
\end{xlist}
is $\cd.$

However, $\Eq$ contains an element $w$
and $2^{\cd}$ elements $y_i\neq w$ such that
\begin{xlist}\item\label{x.joinw}
$y_i\vee y_j=w$ for all distinct $i,j\in I.$
\end{xlist}
\end{proposition}\begin{proof}
The upper bound in the first assertion follows from
the preceding theorem, by the same reasoning
used to get our corollaries to Theorem~\ref{T.gp+md.crit}.
To see without calling on Whitman's Theorem that there does, however,
exist
such a family of cardinality $\cd,$ let $z$ be the discrete equivalence
relation on $\Omega,$ choose an element $p_0\in\Omega$ and a family
of $\cd$ distinct elements $p_i\in\Omega-\{p_0\},$ and for each
$i,$ let $x_i$ be the equivalence relation that relates $p_0$
with $p_i,$ but relates no other distinct elements of $\Omega\<.$

To get the second assertion, let $w$ be the indiscrete
equivalence relation, and let the $y_i$ be all the equivalence
relations on $\Omega$ having exactly two equivalence classes.
\end{proof}

The other way to deal with the fact that the join
operation on $\Eq$ is not continuous in the
topology with subbasis~(\ref{x.Upqc}) is to weaken the topology.
By the discussion following~(\ref{x.Upqc}),
both operations are continuous in the topology having
the sets $U_{p,q}$ as a subbasis of open sets.
This topology is not $\mathrm{T}_1{:}$ for
every $x\in\Eq,$ we see that the closure of $\{x\}$
is the set of all equivalence relations $\leq x.$

As a consequence, the diagonal subset of $\Eq\times\Eq$
is not closed; hence, though arbitrary lattice words
$v$ and $w$ with constants in $\Eq$ still
induce continuous operations on $\Eq,$ it does not follow that the sets
$S_{v=w}=\{x\in\Eq\mid v(x)=w(x)\}$ are closed.

However, let us make
\begin{definition}\label{D.S^leq}
For any lattice $A,$ any set $J,$ any lattice word $v$ in a
$\!J\!$-tuple of variables with constants in $A,$ and
any element $c\in A,$ let
\begin{xlist}\item\label{x.leqset}
$S_{v\leq c}\ =\ \{a=(a_j)_{j\in J}\in A^J\mid
v(a)\leq c\}\ \subseteq\ A^J,$
\end{xlist}
and let us call such sets principal {\em lower} solution sets in $A^J.$
A {\em lower solution set} will mean the intersection of an arbitrary
family of principal lower solution sets.
We shall denote by $L^{\strt[3.3]\leq\mathrm{const}}_{A,\,J}$ the
lattice of all lower solution sets in $A^J.$
\end{definition}

As in Definition~\ref{D.S^=}, we allow the empty intersection, so
that $L^{\strt[3.3]\leq\mathrm{const}}_{A,\,J}$ is
indeed a complete lattice.

Any principal lower solution set $S_{v\leq c}$ in $\Eq^J$ is the
inverse image under the continuous map $v:\Eq^J\rightarrow\Eq$
of the closed set $\{d\mid d\leq c\}\subseteq\Eq,$ hence is
closed, making Lemma~\ref{L.cdB} applicable.
Thus we get

\begin{theorem}\label{T.latleqcrit}
Let $J$ be any set of cardinality $\leq\cd.$
Then $L^{\strt[3.3]\leq\mathrm{const}}_{\Eq,\,J}$ contains
no chains with $>\cd$ jumps.
Hence for any lattice $A$ embeddable in $\Eq,$
$L^{\strt[3.3]\leq\mathrm{const}}_{A,\,J}$ contains
no such chains.\qed
\end{theorem}

To get applications of this result that are not consequences
of Theorem~\ref{T.semilatcrit}, one has to use sets $S_{v\leq c}$
determined by words $v$ that involve both meets and joins;
in fact, that involve meets {\em of} joins.
For if $v$ does not involve meets of joins, we can write
$v=\bigvee_{r=1}^m \bigwedge_{s=1}^{n_r} a_{r,s},$ where
the $a_{r,s}$ are variables and/or constants, and
$m, n_1,\dots,n_m$ are positive integers (some of which may be $1).$
Then we see that the relation $v\leq c$ is equivalent to the
conjunction of the relations
$\bigwedge_{s=1}^{n_r} a_{r,s}\leq c$ $(r=1,\dots,m),$
and each of these can be rewritten
$c\wedge\bigwedge_{s=1}^{n_m} a_{r,s}=\bigwedge_{s=1}^{n_m} a_{r,s}.$
Hence $S_{v\leq c}\in L^{\strt=}_{\Eq_\wedge,\,J},$ and
we are reduced to Theorem~\ref{T.semilatcrit}.
Perhaps workers in lattice theory will be able to
see interesting applications of Theorem~\ref{T.latleqcrit}
that do not reduce to Theorem~\ref{T.semilatcrit}.

We remark that while Theorems~\ref{T.semilatcrit} and~\ref{T.latleqcrit}
show that $L^{\strt=}_{\Eq_\wedge,\,J}$ and
$L^{\strt[3.3]\leq\mathrm{const}}_{\Eq,\,J},$
contain no well-ordered chains of cardinality $>\cd,$
the second paragraph of Proposition~\ref{P.mtvsjn} shows that
$L^{\strt=}_{\Eq,\,J}$ and $L^{\strt=}_{\Eq_\vee,\,J}$ do contain
such chains.

\section{Monoids embeddable in $\Rel.$}\label{S.notRel}

In this section we return to monoids, but focus,
not on $\Se,$ but on its ``wild sibling'' $\Rel,$
the monoid of all binary {\em relations} on
$\Omega,$ under relational composition:
\begin{xlist}\item\label{x.defcomp}
$x\,y\ =\ \{(q,p)\in\Omega\times\Omega\mid(\exists\,r\in\Omega)\ %
(q,r)\in x,\ (r,p)\in y\}.$
\end{xlist}
Note that when specialized to functions, the above definition
of composition, together with our convention that functions act on the
left and are composed accordingly, requires that we
identify each function $f$ with the set of ordered pairs $(f(p),\,p),$
rather than the more usual $(p,\,f(p)).$

Clearly, $\Rel\cong\Rel\op,$ and this monoid contains
$\Se;$ hence it also contains an isomorphic copy of $\Se\op.$
In fact, it is generated by $\Se$ and the natural copy
of $\Se\op;$ for if $R\in\Rel$ is nonempty, then since, as a subset of
$\Omega\times\Omega,$ it has cardinality $\leq\cd,$ we can find
maps $f,g\in\Se$ such that $R=\{(f(p),g(p))\mid p\in\Omega\}.$
Letting $\bar g=\{(p,g(p))\mid p\in\Omega\}\in\Se\op,$
we easily verify that $f\,\bar g=R.$
The empty relation, on the other hand, equals
$\bar f\,g$ for any $f$ and $g$ having disjoint ranges.

The next result shows that monoids embeddable in $\Rel$ are considerably
less restricted than those embeddable $\Se$ or $\Se\op.$
\begin{proposition}\label{P.RelvsSe}
Let $w$ denote
the indiscrete equivalence relation on $\Omega,$ and
\textup{(}as in the proof of Proposition~\ref{P.mtvsjn}\textup{)}
let $y_i$ $(i\in I)$ be the $2^{\cd}$ distinct equivalence relations on
$\Omega$ having exactly two equivalence classes.
Then as members of $\Rel,$ these satisfy
\begin{xlist}\item\label{x.yiyi}
$\!y_i\<y_i = y_i$ $(i\in I),$
\end{xlist}
\begin{xlist}\item\label{x.yiyjyi}
$y_i\<y_j\<y_i = w$ $(i\neq j$ in $I).$
\end{xlist}
Hence $L^{\strt=}_{\Rel, 1}$ contains a well-ordered chain of
cardinality $2^{\cd}.$
\end{proposition}\begin{proof}
(\ref{x.yiyi}) is immediate.
The verification of~(\ref{x.yiyjyi}) is routine, but the reader may
find the following way of visualizing it helpful.
Given $y_i$ and $y_j,$ picture a Venn diagram for $\Omega,$
divided by a vertical line representing the partition
into the two equivalence classes of $y_i$ and a horizontal
line representing the partition into the equivalence classes of $y_j.$
Given $(p,q)\in w=\Omega\times\Omega,$ which we wish to show
belongs to $y_i\<y_j\<y_i,$ we may assume by adjusting
our diagram that $p$ lies in the upper left-hand box.
Since $y_i$ and $y_j$ each have two equivalence classes, at least
one of the lower boxes and at least one of the right-hand boxes
are nonempty, and since $y_i\neq y_j,$ the lower right-hand box
is not the only nonempty box other than the upper left-hand one.
It is now easy to see that wherever $q$ may lie in our diagram,
we can get from $p$ to it by crossing the vertical line at most
once and the horizontal line at most once, hence that $(p,q)$
lies in either $y_i\<y_j$ or $y_j\<y_i.$
In either case, it lies in $y_i\<y_j\<y_i,$ as claimed.

However,~(\ref{x.yiyi})
shows that $y_i\<y_i\<y_i\neq w,$ and the contrast between this
inequality and~(\ref{x.yiyjyi}) allows us to get, as in the proofs of
Corollaries~\ref{C.cntrlzr} and~\ref{C.235}, a well-ordered chain
of cardinality $2^{\cd}$ in $L^{\strt=}_{\Rel, 1.}$
\end{proof}

The composition operation~(\ref{x.defcomp})
of $\Rel$ resembles the join operation of $\Eq$
in being continuous in the topology on $\Rel$ with subbasis
of open sets consisting of the sets $U_{p,q}$ (defined
in~(\ref{x.Upqc})), but not in the Hausdorff topology with subbasis
of open sets $U_{p,q}$ and $^c U_{p,q}.$
Indeed, though composition of relations is simpler to
describe than the join of
equivalence relations, it is still true that no specification of whether
some {\em finite} number of pairs belong to each of two
relations can tell us that $(p,q)$ does not belong to their composite.
For $\Rel,$ there is no analog of our ploy of restricting attention
to the meet operation of $\Eq;$ so
let us go directly to the weaker topology.
As in the case of $\Eq,$ we can conclude that

\begin{lemma}\label{L.relsol}
For $J$ any set of cardinality $\leq\cd,$
$L^{\strt[3.3]\leq\mathrm{const}}_{\Rel,\,J}$ contains no chains
with $>\cd$ jumps.\qed
\end{lemma}

Unfortunately, I don't see how to use this result in studying
the structures of
monoids embeddable in the monoid $\Rel,$ for the definition
of $L^{\strt[3.3]\leq\mathrm{const}}_{\Rel,\,J}$ makes use of
the ordering of $\Rel$ by inclusion, and,
unlike the order relation on a lattice, this
has no description in terms of its algebra operations.
Of course, if we throw in order as additional structure,
then Lemma~\ref{L.relsol} yields restrictions on embeddability.
Let us record this and another consequence
of Lemma~\ref{L.relsol}, without looking for examples,
then return to the question of embeddability as pure monoids.

By a {\em partially ordered monoid} $(M,\leq),$ let us understand a
monoid $M$ given with a partial ordering $\leq,$ such that for
all $x,\,y,\,z\in M,$
\begin{xlist}\item\label{x.ordmonoid}
$\!x\leq y\implies xz\leq yz,\\[2pt]
\!x\leq y\implies zx\leq zy.$
\end{xlist}
An {\em embedding} of partially ordered monoids
$f:(M,\leq_M)\rightarrow (N,\leq_N)$ will mean a monoid embedding
$f:M\rightarrow N$ such that the partial ordering on $M$
induced by $\leq_N$ under $f$ is precisely $\leq_M.$
Clearly, Lemma~\ref{L.relsol} gives
\begin{corollary}\label{C.ordmonoid}
If a partially ordered monoid $(M,\leq)$ is embeddable in
$(\Rel,\subseteq),$ then for any set $J$ of cardinality $\leq\cd,$
the complete lattice $L^{\strt[3.3]\leq\mathrm{const}}_{(M,\leq),\,J}$
contains no chains with $>\cd$ jumps.\qed
\end{corollary}

To formulate another consequence of Lemma~\ref{L.relsol},
let us define, for $M$ a monoid
with a distinguished element $z,$ $J$ a set, and $v$ any monoid
word in a $\!J\!$-tuple of variables and arbitrary
constants from $M,$ the set
$S_{v=z}=\{a=(a_j)_{j\in J}\in M^J\mid v(a)=z\}\subseteq M^J,$
and let $L^{\strt[3.3]=z}_{(M,z),\,J}\subseteq \Pw[M^J]$ denote the
lattice of all intersections of families of sets of this sort.
If $(M,z)$ and $(M',z')$ are such pairs,
a homomorphism $(M,z)\rightarrow (M',z')$ will mean
a monoid homomorphism $M\rightarrow M'$ carrying $z$ to $z'.$

\begin{corollary}\label{C.empty}
Let $M$ be a monoid with a zero element $z$ \textup{(}an element
satisfying $zx=z=xz$ for all $x\in M),$ such that
$(M,z)$ is embeddable in $(\Rel,\emptyset),$
and let $J$ be a set of cardinality $\leq\cd.$
Then the complete lattice $L^{\strt[3.3]=z}_{(M,z),\,J}$
contains no chains with $>\cd$ jumps.\qed
\end{corollary}\begin{proof}
In $\Rel^J,$ any set of the form $S_{v=\emptyset}$ can clearly
also be described as $S_{v\leq\emptyset},$ whence the assertion
follows immediately from Lemma~\ref{L.relsol}.
\end{proof}

A way of getting restrictions on embeddability in $\Rel$ without
bringing in additional structure is to note, as we did
in~\S\ref{S.notSym} for
$\Se,$ that any monoid homomorphism from a group into $\Rel$ will land
in the group of invertible relations, which is again $\Sym.$
Hence our restrictions on groups embeddable
in $\Sym$ are also restrictions on the groups of invertible
elements of monoids embeddable in $\Rel.$

But in fact, we can use the group $\Sym$ of invertible elements
of $\Rel$ in a way that brings in noninvertible elements as well.
Note that the function topology on $\Se,$ and hence
on $\Sym,$ is the restriction thereto both of the topology
of~(\ref{x.Upqc}), and of the weaker
topology having only the sets $U_{p,q}$ as subbasis of open sets.
(In $\Se,$ the set $^c U_{p,q}$ can be written
as $\bigcup_{p'\neq p} U_{p',q},$ so the topology generated by
the sets $U_{p,q}$ also contains the complementary sets.)
This immediately gives case~(i) of the next lemma.
Case~(ii), the one we shall make use of, is more surprising.
(We shall not use the final parenthetical strengthening of~(ii).)

\begin{lemma}\label{L.sym+rel}
The restrictions of the monoid multiplication of $\Rel$ to maps
\begin{xlist}\item\label{x.sym+rel}
$\Sym\times\Rel\rightarrow\Rel$ \quad and
\quad $\Rel\times\Sym\rightarrow\Rel$
\end{xlist}
are continuous if we put the function topology on $\Sym,$ and
put on $\Rel$ the topology with subbasis of open sets consisting
either of\\
\textup{(i)} the sets $U_{p,q},$ or\\
\textup{(ii)} the sets $U_{p,q}$ and $^c U_{p,q}.$
\textup{(}This remains true if we replace $\Sym$
by $\Se$ in the second map of~\textup{(\ref{x.sym+rel})},
though not in the first.\textup{)}
\end{lemma}\begin{proof}
As noted, continuity in the topology determined
by~(i) follows from the continuity of
the multiplication of $\Rel$ in that topology.
This also gives half of continuity in the topology specified in~(ii),
namely openness of the inverse images of the sets $U_{p,q}.$
We shall prove the corresponding statement for $^c U_{p,q}$ for the
second map in~(\ref{x.sym+rel}) with $\Sym$ replaced by $\Se.$
The case where $\Sym$ is left unchanged
follows immediately, and the first statement of~(\ref{x.sym+rel})
then follows by reversing coordinates in ordered pairs.

Note that a necessary and sufficient condition for a
composite $x\<a$ $(x\in\Rel,\,a\in\Se)$ to lie in $^c U_{p,q},$
i.e., not to contain $(p,q),$ is that for the unique $r\in\Omega$
such that $(r,q)\in a,$ we have $(p,r)\notin x.$
Hence the inverse image of $^c U_{p,q}$ under the above composition
map is the union
$\bigcup_{r\in\Omega}\,({^c U_{p,r}}\times (U_{r,q}\cap\Se)).$
This set is open in $\Rel\times\Se,$ as claimed.

To get the negative part of the final parenthetical assertion,
assume without loss of generality that $\omega\subseteq\Omega\<.$
For each $n\in\omega,$ let $x_n\in\Rel$ be the partial
function which sends $0$ to $n$ and does nothing else, and let
$a\in\Se$ be the function sending all elements to $0.$
Then for all $n,$ $a\,x_n=x_0,$
hence $\lim_{n\rightarrow\infty} a\,x_n=x_0.$
On the other hand, $\lim_{n\rightarrow\infty} x_n=\emptyset,$ the empty
relation, which when left multiplied by $a$ gives $\emptyset\neq x_0.$
So composition is not continuous in the Hausdorff topology of~(ii).
\end{proof}

From the operations~(\ref{x.sym+rel}) and the group operations of
$\Sym,$ we can form words in a mixture of $\!\Sym\!$-
and $\!\Rel\!$-valued elements; but
any such word can involve, at most, {\em either} one occurrence
of a $\!\Rel\!$-valued variable
{\em or} one occurrence of a non-$\!\Sym\!$-valued constant, since the
operations~(\ref{x.sym+rel}) do not allow the multiplication
of two non-$\!\Sym\!$-valued elements.
(Here by a ``$\!\Rel\!$-valued variable'' I mean a variable that
is allowed to range over all of $\Rel,$ taking on both invertible
and noninvertible values.
Constants, on the other hand, have specific values rather than ranges,
so for these the relevant concept is that of a non-$\!\Sym\!$-valued
constant, i.e., a noninvertible element of $\Rel.)$
By Lemma~\ref{L.sym+rel}, if we take a family
of such words, possibly involving many $\!\Rel\!$-valued variables
and non-$\!\Sym\!$-valued constants altogether (though with at most
one of these per word), and a set of equations in these words,
then its solution set is closed in the Hausdorff topology
defined by the subbasis~(ii).
Hence we can apply Lemma~\ref{L.cdB} and get restrictions on the
lattice of such solution sets, which imply the same restrictions
on the corresponding lattice obtained from any monoid
embeddable in $\Rel.$

It is not clear to me whether so allowing more than one
$\!\Rel\!$-valued variable or constant in our system of
equations actually contributes to the generality of this result.
Note that if we have a system of equations in several such variables
and constants, then any equation involving different $\!\Rel\!$-valued
variables on the two sides, say $x$ on the left and $y$ on the right,
will allow us to solve for $y$ in terms of $x$ and
the $\!\Sym\!$-valued variables and constants (since the latter can
all be inverted and brought to the left side of the equation); and we
can then substitute the resulting expression for all occurrences
of $y$ in the remaining equations, and so eliminate $y$ from the system.
On the other hand, our interest is not in the solution set of a
{\em single}
system of equations but in the relation between solution sets of
many such systems, and a variable that can be eliminated
from one of these will not in general be eliminable from all of them.
So if the consequences of our conditions on lattices of
solution sets can indeed be reduced to the case where there
is only one $\!\Rel\!$-valued variable, the argument by
which this reduction is done may be nontrivial.

Leaving it to others to determine whether such a reduction is
possible, I will, for simplicity, record here only the statements for
systems with at most one $\!\Rel\!$-valued-variable.

\begin{theorem}\label{T.Relcrit}
Let $J$ be any set of cardinality $\leq\cd.$

Let $L$ denote the lattice of solution sets in
$\Sym^J\times\Rel$ of systems of equations in a $\!J\!$-tuple of
$\!\Sym\!$-valued variables, and a single $\!\Rel\!$-valued variable
which appears at most once on each side of any equation, together
with arbitrary $\!\Sym\!$-valued constants.

Let $L'$ similarly denote the lattice of solution sets in
$\Sym^J$ of systems of equations having on each
side a word in a $\!J\!$-tuple of $\!\Sym\!$-valued variables,
arbitrary $\!\Sym\!$-valued constants,
and at most one occurrence of a $\!\Rel\!$-valued constant.

Then neither of these lattices contains a chain with $>\cd$ jumps.
Hence for any monoid $M$ embeddable in $\Rel,$
the corresponding restrictions hold, with the group $U(M)$ of invertible
elements of $M$ taking the place of $\Sym.$
\qed
\end{theorem}

Simple examples of equations of the sort arising above are
stabilizer relations $gx=x$ and $xg=x$ $(x\in M,\,g\in U(M)),$ and more
generally, $gxg'=x$ $(x\in M,\,g,g'\in U(M)).$
Here is an example based on relations $gx=x.$

Recall that a {\em left zero} element in a monoid $M$
means an element $z$ satisfying $zx=z$ for all $x\in M.$

\begin{corollary}\label{C.gzz'}
Let $I$ be a set,
let $G$ be the restricted direct product of an $\!I\!$-tuple
of copies of the group $Z_2,$ with generators $g_i$ $(i\in I)$
\textup{(}i.e., the additive group of the vector space with
basis $\{g_i\mid i\in I\}$ over the field of two elements\textup{)},
and let $M$ be the monoid whose group of invertible elements is
this group $G,$ and whose other elements, denoted
$z_i,\,z_i'$ $(i\in I),$ are left zero elements whose
behavior under the left action of $G$ is described by
\begin{xlist}\item\label{x.gzz'}
Left multiplication by $g_i$ interchanges $z_i$ and $z'_i,$
and fixes all $z_j$ and $z'_j$ with $j\neq i.$
\end{xlist}

Then $M$ is embeddable in $\Rel$ if and only if $\cd[I]\leq\cd.$
\end{corollary}\begin{proof}
It is straightforward to verify that the above operations define
a monoid, which has cardinality $\leq\max(\cd[I],\aleph_0),$ and so is
embeddable in $\Se\subseteq\Rel$ if $\cd[I]\leq\cd.$
To prove the ``only if'' part of the conclusion, assume
$I$ is a cardinal $\kappa>\cd.$
For each $\alpha\in\kappa,$ let $S_\alpha=
\{x\in M\mid(\forall\,\beta>\alpha)\ \linebreak[1]g_\beta\,x=x\}.$
Note that $z_\beta\in S_\alpha$ if and only
if $\beta\leq\alpha,$ so the $S_\alpha$ form a well-ordered
chain of cardinality $\kappa>\cd.$
By the ``$\!L\!$'' case of Theorem~\ref{T.Relcrit} (with $J=0),$
this precludes embeddability in $\Rel.$
\end{proof}

In proving the above corollary we could, alternatively,
have used the $z_\alpha$ as monoid-valued constants,
and let the $g$ in the relations $g\,z_\alpha=z_\alpha$ be an
invertible-element-valued variable, getting a system to which
the ``$\!L'\!$'' case of Theorem~\ref{T.Relcrit} with $J=1$ applied.

\begin{question}\label{Q.Rel}
Can one give stronger necessary conditions for
embeddability of a monoid $M$ of cardinality $\leq 2^{\cd}$
in $\Rel$ than those of Theorem~\ref{T.Relcrit}?

In particular, suppose we define the left, right, and
$\!2\!$-sided stabilizers of an element $x$ of a monoid $M$
as $\{y\in M\mid yx=x\},$ $\{y\in M\mid xy=x\},$ and
$\{(y,y')\in M^2\mid yxy'=x\}.$
Does $\Rel$ have any chains of intersections of such stabilizers
with $>\cd$ jumps?
What if we restrict $y$ and $y'$ here to {\em left
invertible} or to {\em right invertible} elements?

Are there any monoids of cardinality $\leq 2^{\cd}$
having no invertible elements other than $1$
\textup{(}or having $\leq\cd$ invertible
elements\textup{)} which are not embeddable in $\Rel$?

For $\cd[I]>\cd,$ is the monoid of Corollary~\ref{C.gzz'}
ever embeddable {\em as a semigroup} in $\Rel$?
\end{question}

One kind of noninvertible elements one might look at in approaching
these questions are the {\em idempotents}, since for $e$ idempotent,
the solution sets of $xe=x$ and $ex=x$ are particularly natural objects.
Idempotents in $\Rel$ come in more forms than one might expect.
Obvious examples are equivalence relations,
subsets of the identity relation, graphs of retractions
of $\Omega$ onto subsets, and the opposites of such graphs;
these four constructions can also be mixed in fairly natural ways.
For less obvious examples, note that (i)~for any
partial ordering, $\leq,$ on $\Omega$ with no jumps (e.g.,
the ordinary ordering on the set of rational numbers), the graph
of the relation ``$\!<\!$'' is idempotent, though it has trivial
intersection with the identity relation; and (ii)~for any nondisjoint
subsets $X$ and $Y$ of $\Omega,$ the set
$X\times Y\subseteq\Omega\times\Omega$ is an idempotent relation.

We have not written down the positive embeddability result
analogous to Theorem~\ref{T.funct}, i.e., the
embeddability in $\Rel$ of monoids $F(2^{\cd})$ for appropriate
functors $F,$ since this follows from Theorem~\ref{T.funct} and
the inclusion $\Se\subseteq\Rel.$
It would, of course, be of interest if one could get {\em stronger}
results of this sort for $\Rel$ than for $\Se.$
I also don't know the answer to

\begin{question}\label{Q.doubleRel}
Is $\Rel\cP\Rel$ embeddable as a monoid in $\Rel$?

\textup{(}If so, then by the usual argument, $\Rel$ in fact contains a
coproduct of $2^{\cd}$ copies of itself.\textup{)}
\end{question}

Here are some partial positive results on that question.

Let $1_\Omega\in\Rel$ denote the identity element,
that is, the diagonal subset of $\Omega\times\Omega\<.$
Let us say that elements $g,\,h\in\Rel$ ``differ off the
diagonal'' if $g-1_\Omega\neq h-1_\Omega$ (where ``$\!-\!$'' denotes
set-theoretic difference); and for a subset $X\subseteq\Rel,$
let us say ``the members of $X$ are
distinguishable off the diagonal'' if every pair of distinct
elements of $X$ differs off the diagonal.

\begin{lemma}\label{L.doubleRel}
There exists a monoid homomorphism $f:\Rel\cP\Rel\rightarrow\Rel$ with
the property that whenever $M$ and $N$ are submonoids of $\Rel,$ each
of which has the property that its members are distinguishable off the
diagonal, then the restriction
of $f$ to a homomorphism $M\cP N\rightarrow\Rel$ is an embedding.
\end{lemma}\begin{proof}
Let us understand a ``relational action'' of a monoid $M$ on a set $X$
to mean a homomorphism $M\rightarrow\Rel[(X)].$
Paralleling the proof of Lemma~\ref{L.double}(ii), we shall construct
a relational action of $\Rel\cP\Rel$ on the disjoint union of a
family of copies of $\Omega\times\omega,$
where the family is again indexed by the group
$\Sym[_{\mathrm{fin}}(\Omega\times\omega)],$
such that the restriction of this action to any submonoid
$M\cP N$ as in the statement of the lemma is faithful.

Let the ``natural relational action''
of $\Rel$ on $\Omega\times\omega$ be defined to take each $g\in\Rel$
to the relation
on $\Omega\times\omega$ consisting of all pairs $((q,k),\,(p,k))$
with $(q,p)\in g$ and $k\in\omega.$
Let $\alpha,\,\beta$ be the two coprojections
$\Rel\rightarrow\Rel\cP\Rel.$
For each $g\in\Rel,$ let us send $\alpha(g)$ to the natural
relational action of $g$ on each copy
of $\Omega\times\omega,$ while letting $\beta(g)$ act
on the copy of $\Omega\times\omega$ indexed by each
$t\in\Sym[_{\mathrm{fin}}(\Omega\times\omega)]$ via
the conjugated relation $t\,g\,t^{-1}.$
To complete our proof, we need to show that for
$M$ and $N$ as in the statement of the lemma,
if~(\ref{x.nmfmg}) and~(\ref{x.nmfmh}) are distinct elements of
$M\cP N\subseteq\Rel\cP\Rel,$ then
there exists $t\in\Sym[_{\mathrm{fin}}(\Omega\times\omega)]$
such that the two relations~(\ref{x.2wt+-})
on $\Omega\times\omega$ are distinct.

As in our previous argument, we assume that the length
of~(\ref{x.nmfmg}), which we will call $n,$ is at least the length
of~(\ref{x.nmfmh}).
Moreover, we assume that if these lengths are equal, and if
$\!\alpha\!$s and $\!\beta\!$s occur in the same positions in both
expressions, then at the first position $k$ from the right
where~(\ref{x.nmfmg}) and~(\ref{x.nmfmh}) differ, $g_k$ contains
some {\em nondiagonal} ordered pair which does not lie in $h_k.$
Since $g_k$ and $h_k$ differ off the diagonal, we can achieve this by
interchanging~(\ref{x.nmfmg}) and~(\ref{x.nmfmh}) if necessary.
We now choose, for each $k\in\{1,\dots,n\},$ a nondiagonal
pair $(q_k,p_k)\in g_k$ (which must exist, because
$g_k$ is not the identity, hence is distinguishable from the
identity off the diagonal), using, when possible, a pair not
also contained in $h_k;$ and we define $t$
as in~(\ref{x.transpose}), except that wherever~(\ref{x.transpose})
shows an element $g_k(p_k),$ we now use $q_k.$
The same considerations as in the proof of Lemma~\ref{L.double}
show that the first product
in~(\ref{x.2woAt}), but not the second, contains the
pair $((q_n,\,n{+}1),\ (p_0,0));$ so these two products are distinct
relations on $\Omega\times\omega,$ completing the proof.
\end{proof}

Here are some consequences.

\begin{corollary}\label{C.doublefinetc}
\textup{(i)}~ $\Se\cP\Se\op$ is embeddable in $\Rel.$

\vspace{2pt}
\textup{(ii)}~ Let $\Rel_{\geq 1}\subseteq\Rel$ denote the
submonoid of relations that contain the diagonal
\textup{(}the reflexive relations\textup{)}.
Then $\Rel_{\geq 1}\,\cP\,\Rel_{\geq 1}$ is embeddable
in $\Rel_{\geq 1}.$
Hence the coproduct of $2^{\cd}$ copies
of $\Rel_{\geq 1}$ is embeddable in $\Rel_{\geq 1}.$

\vspace{2pt}
\textup{(iii)}~ Let $\Rel[^{\mathrm{semi}}(\Omega)]$
denote the underlying semigroup of $\Rel.$
Then the {\em semigroup} coproduct\linebreak[4]
$\Rel[^{\mathrm{semi}}(\Omega)]\cP\Rel[^{\mathrm{semi}}(\Omega)]$
is embeddable in $\Rel[^{\mathrm{semi}}(\Omega)].$
Hence the semigroup coproduct of $2^{\cd}$ copies
of $\Rel[^{\mathrm{semi}}(\Omega)]$
is embeddable in $\Rel[^{\mathrm{semi}}(\Omega)].$
\end{corollary}\begin{proof}
(i) is immediate from the lemma, since both $\Se$ and the natural
copy of $\Se\op$ in $\Rel$ (gotten by taking the opposite relations to
all members of $\Se)$ have the property that their members are
distinguishable off the diagonal.

In~(ii), Lemma~\ref{L.doubleRel} immediately gives
embeddability of $\Rel_{\geq 1}\cP\Rel_{\geq 1}$ in $\Rel.$
Moreover, we see that the construction of that lemma
takes reflexive relations to
reflexive relations, so the embedding lands in $\Rel_{\geq 1}.$
The method of proof of Theorem~\ref{T.bigcP} now allows one to work
one's way up to the coproduct of $2^{\cd}$ copies.

To prove~(iii) let us begin by constructing an embedding
of $\Rel[^{\mathrm{semi}}(\Omega)]$ in the semigroup of relations
on a set of the same cardinality
as $\Omega,$ by relations that are distinguishable from
each other and from the identity relation off the diagonal.

We first have to declaw the empty relation; so let $z$ be an element
not in $\Omega,$ and let us embed $\Rel$ in $\Rel[(\Omega\cup\{z\})]$
by sending each relation $g$ to $g\cup\{(z,z)\};$ this is an
embedding of monoids whose image consists of nonempty relations.

We now map $\Rel[^{\mathrm{semi}}(\Omega\cup\{z\})]$ into
$\Rel[^{\mathrm{semi}}((\Omega\cup\{z\})\times 2)]$ by
\begin{xlist}\item\label{x.times2}
$g\ \mapsto\ \{((q,i),\,(p,j))\mid (q,p)\in g;\ i,j\in 2\}.$
\end{xlist}
This construction is easily seen to respect composition,
and to take distinct nonempty relations to relations that differ
off the diagonal both from each other and from $1_{\Omega\cup\{z\}}.$
The image of the composite map
$\Rel[^{\mathrm{semi}}(\Omega)]\rightarrow
\Rel[^{\mathrm{semi}}(\Omega\cup\{z\})]\rightarrow
\Rel[^{\mathrm{semi}}((\Omega\cup\{z\})\times 2)]$
is in particular a subsemigroup
$S\subseteq\Rel[^{\mathrm{semi}}((\Omega\cup\{z\})\times 2)]$
isomorphic to $\Rel[^{\mathrm{semi}}(\Omega)].$
Applying Lemma~\ref{L.doubleRel} with $S\cup\{1\}$
in the role of both $M$ and $N,$ we get an embedding of monoids
$(S\cup\{1\})\linebreak[1]\cP\linebreak[3](S\cup\{1\})\rightarrow\Rel,$
which, restricted to the subsemigroup generated by the two
copies of $S,$ gives the desired embedding of semigroups.
As before, the method of \S\ref{S.jiggle} allows us to push this
up to an embedding of a $\!2^{\cd}\!$-fold coproduct of copies of
$\Rel[^{\mathrm{semi}}(\Omega)]$ in $\Rel[^{\mathrm{semi}}(\Omega)].$
\end{proof}

It would be nice if we could carry the idea of part~(iii) of the above
corollary further.
If we could find a {\em monoid} embedding
$\varphi: \Rel\rightarrow\Rel[(\Omega')]$ with $\cd[\Omega']=\cd,$
such that images of distinct elements were distinguishable
off the diagonal, then an application of Lemma~\ref{L.doubleRel}
would give us our desired embedding $\Rel\cP\Rel\rightarrow\Rel.$
But when we attempt to construct such a $\varphi,$ we run into
difficulty trying to simultaneously

(a)~make it carry the identity relation $1_\Omega$ to $1_{\Omega'},$

(b)~handle relations which are properly contained in $1_\Omega,$ and

(c)~handle relations
which are both infinitely-many-to-one and one-to-infinitely-many.\\
Part~(iii) of the above corollary showed that we could get
an embedding if we dropped the requirement~(a).
The following result (strengthening part~(ii) of that corollary)
shows the same if we instead drop~(b), and the result after that
will do the same for~(c).

\begin{lemma}\label{L.offdiag}
There exists a monoid homomorphism $f':\Rel\cP\Rel\rightarrow\Rel$ with
the property that whenever $M$ and $N$ are submonoids of $\Rel,$
neither of which contains a proper subrelation of $1_\Omega,$
then the restriction
of $f'$ to a homomorphism $M\cP N\rightarrow\Rel$ is an embedding.
\end{lemma}\begin{proof}
Let $\varphi:\Rel\rightarrow\Rel[(\Omega^2)]$ take each
$g\in\Rel$ to $\{((p,p'),(q,q'))\mid (p,q),(p',q')\in g\}.$
It is straightforward to verify that this is a monoid homomorphism.
We shall show that if $g,h\in\Rel$ are distinct elements which
are not both subrelations of $1_\Omega,$ then $\varphi(g)$
and $\varphi(h)$ are distinguishable off the diagonal.
This, together with Lemma~\ref{L.doubleRel}, yields the desired result.

First, suppose $g$ and $h$ are themselves distinguishable off
the diagonal.
Then without loss of generality we may assume that $(p,q)$ is
contained in $g$ but not $h$ for some $p\neq q,$ and we see that
$((p,p),(q,q))$ belongs to $\varphi(g)$ but not $\varphi(h).$
On the other hand, if $g$ and $h$ are not distinguishable off
the diagonal, then since they are not both subrelations of
$1_\Omega,$ there must be
some $(p,q)$ with $p\neq q$ that is contained in both of them.
Also, since $g\neq h,$ some $(r,r)$ $(r\in\Omega)$ will belong
to one of them but not the other; say to $g$ but not $h.$
Then we see that $((p,r),(q,r))$ is a nondiagonal
element belonging to $\varphi(g)$ but not $\varphi(h).$
\end{proof}

Finally, let us see what we can get if we sacrifice~(c).

Given $g\in\Rel$ and $X\subseteq\Omega,$ let us define the ``image'' set
\begin{xlist}\item\label{x.gX}
$g\<X\ =\ \{q\in\Omega\mid(\exists\,p\in X)\ (q,p)\in g\},$
\end{xlist}
and let
\begin{xlist}\item\label{x.Relleftfin}
$\Rel[_{\mathrm{fin}\leftarrow 1}(\Omega)]\ =\ %
\{g\in\Rel\mid(\forall\,p\in\Omega)\ \ g\{p\}$ is finite$\!\}.$
\end{xlist}
The arrow points to the left to show that we are defining this
set in terms of the left action~(\ref{x.gX}).
We can clearly
also describe~(\ref{x.Relleftfin}) as the set of $g$ such that
for every finite $X\subseteq\Omega,$ the set $g\<X$ is again finite.
This shows $\Rel[_{\mathrm{fin}\leftarrow 1}(\Omega)]$
to be a submonoid of $\Rel,$ and it clearly
acts faithfully -- by functions, not relations -- on $\Pf.$
This gives us an embedding $\Rel[_{\mathrm{fin}\leftarrow 1}(\Omega)]
\rightarrow\Se[(\Pf)]\cong\Se.$
Of course, we also have
$\Se\subseteq\Rel[_{\mathrm{fin}\leftarrow 1}(\Omega)];$ so
embeddability of a monoid in $\Rel[_{\mathrm{fin}\leftarrow 1}(\Omega)]$
and in $\Se$ are equivalent.
This is the first assertion of the next lemma, and by
the results of \S\ref{S.jiggle} it implies the second.

\begin{lemma}\label{L.RlfSe}
$\Rel[_{\mathrm{fin}\leftarrow 1}(\Omega)]$ and $\Se$
are each embeddable in the other.

Hence $\Rel[_{\mathrm{fin}\leftarrow 1}(\Omega)]$ contains a coproduct
of $2^{\cd}$ copies of itself as a monoid.\qed
\end{lemma}

Let me record here a
curious construction which I thought, at one point, would give
a more elegant proof of our semigroup-embedding
result, Corollary~\ref{C.doublefinetc}(iii).
This did not quite work; but perhaps it is nonetheless of interest.

Consider the map $\theta:\Rel[^{\mathrm{semi}}(\Omega)]\rightarrow
\Rel[^{\mathrm{semi}}(\Pf)]$ defined by
\begin{xlist}\item\label{x.pfim}
$\theta(g)\ =\ \{(t,s)\in(\Pf)^2 \mid t\subseteq gs\}.$
\end{xlist}
It is easy to verify that $\theta$ is a semigroup homomorphism.
Moreover, every $\theta(g)$ contains all pairs
$(\emptyset,t)$ with $t\in\Pf,$
and hence differs from $1_{\Pf}$ off the diagonal.
For most pairs of distinct $g,h\in\Rel,$ one finds
that $\theta(g)$ and $\theta(h)$ differ in the nondiagonal pairs of
one of the forms $(\{p\},\{q\}),$
$(\{p\},\{p,q\}),$ or $(\{p,q\},\{p\})$ that they contain.
There are exceptions, however; for instance, if
$p\in\Omega$ and we take $g=\{p\}\times\Omega,$
$h=\{p\}\times(\Omega-\{p\});$ then $\theta(g)$ and $\theta(h)$
differ only with regard to the diagonal pair $(\{p\},\{p\}).$
Hence this construction is not itself a substitute for
the one used in proving Corollary~\ref{C.doublefinetc}(iii).
We can cure the problem using the
``$\!\Rel\rightarrow\Rel[(\Omega\cup\{z\})]\!$''
trick that was used (to cure a different problem) in the existing proof
of that statement (and after this is done, the images of any two
relations in fact differ in the pairs $(\{p\},\{q,z\})$ that they
contain); but I wouldn't call the resulting proof more elegant than the
one we used.

We end this section with
a different sort of self-embeddability question for $\Rel.$
\begin{question}\label{Q.eRe}
If $e\in\Rel$ is an idempotent, and we regard
$eRe=\{x\in\Rel\mid exe=x\}$
as a monoid with identity element $e,$ is this monoid always
embeddable in $\Rel$?
\textup{(}Cf.\ remarks on idempotents following
Question~\ref{Q.Rel}.\textup{)}
\end{question}

\section{Restrictions on $\!K\!$-algebras embeddable in
$\En.$}\label{S.notEn}

Recall that $V$ denotes a vector space with basis $\Omega$
over a field $K.$
Since the endomorphism algebra $\En$ is a $\!K\!$-linear analog of
the monoid $\Se,$ we would hope to get restrictions on
associative unital $\!K\!$-algebras embeddable in $\En$
parallel to our restrictions on monoids embeddable in $\Se.$
We can do this -- except that, where we would like to bound
the number of jumps in a chain of solution sets by $\cd=\dim_K(V),$
I only know how to bound it by $\cd[V]=\max(\cd,\cd[K]).$
The following theorem is proved exactly like Theorem~\ref{T.gp+md.crit},
using the fact that addition and composition of members of
$\En,$ and multiplication of these maps by members of $K,$ are
continuous in the function topology on $\En,$ regarded as
a subset of $\Se[(V)].$

\begin{theorem}\label{T.K-alg.crit}
Let $J$ be any set of cardinality $\leq\cd[V].$
Then $L^{\strt=}_{\En,\,J}$ contains no chains with $>\cd[V]$ jumps.
Hence for any $\!K\!$-algebra $A$ embeddable in $\En,$ the lattice
$L^{\strt=}_{A,\,J}$ contains no such chains; in
particular, it contains no well-ordered or reverse-well-ordered
chains of cardinality $>\cd[V].$\qed
\end{theorem}

In the hope of reducing the bound $\cd[V]$ to $\cd,$ we might try
replacing the function topology on $\En$ by some topology with
a smaller basis of open sets; say one that defines its subbasic open
sets not by considering the values of elements of $\En$ at arbitrary
elements of $V,$ but only at the elements of our basis $\Omega\<.$
Unfortunately, each of these images still has $\cd[V]$ possible values.
However, a linear restriction on these images corresponds to
a proper subspace of the space $V$ of possible values, suggesting
that the vector space dimension should still bound lengths of chains.
On the other hand, the conditions on the coordinates of our elements
induced by ring-theoretic relations are not necessarily linear.
Perhaps one should seek bounds on the lengths of chains of
solution sets by methods of algebraic geometry.
Or perhaps one can get stronger results for relations that are
multilinear, such as centralizer and annihilator relations, than
for general relations.

In another direction, if $K$ is a topological field such as the real
or complex numbers, perhaps we could use the topology of that
field instead of the discrete topology on our coordinates, and
replace $\cd[V]$ in the above theorem by
$\max(\cd,\kappa)$ where $\kappa$ is the least
cardinality of a basis for the topology of $K.$
We record these problems as

\begin{question}\label{Q.K-dim}
In Theorem~\ref{T.K-alg.crit}, can the bound $\cd[V]$ be
lowered to $\cd$?
If not in general, what if we restrict attention to the lattice
determined by $\!K\!$-linear or $\!K\!$-affine relations?
Can one at least improve the bound $\cd[V]$ if $K$ admits
a structure of topological field with a basis of $<2^{\cd}$ open sets?
\end{question}

For one sort of system of relations,
we can indeed get the expected bound.

\begin{lemma}\label{L.orthog}
Any set of nonzero {\em pairwise orthogonal idempotent}
elements of $\En$ has cardinality $\leq\cd\<.$
Hence the same is true in any $\!K\!$-algebra embeddable in $\En.$

In particular, $\En$ does not contain a direct product
of $>\cd$ nontrivial $\!K\!$-algebras.
\end{lemma}\begin{proof}
Given an infinite family $S$ of nonzero pairwise orthogonal
idempotents, the images of these as endomaps of $V$ form
a set of subspaces of $V$ whose sum is direct.
Hence the dimension of that sum is at least $\cd[S];$
so $\cd[S]\leq\mathrm{dim}(V)=\Omega\<.$

To see the final sentence, note that in a direct product of
$\!K\!$-algebras $\prod_{i\in I} A_i,$ if $e_i$ denotes the element with
$\!i\!$-component $1$ and all other components $0,$ then these
are pairwise orthogonal idempotents.
\end{proof}

On the other hand, some ways in which the behavior of $\En$ for $K$
large is indeed different from that of
$\Se$ will be noted at the end of~\S\ref{S.jiggleEn}.

Wehrung, paralleling his result mentioned earlier
that $\Se$ cannot be embedded in $\Se\op,$ also shows in~\cite{FW_op}
that that $\En$ cannot be embedded in $\En\op.$
In fact, he shows that $\Se$ cannot be embedded in the
underlying multiplicative semigroup of $\En\op,$
yielding both results!

\section{Other directions for generalization.}\label{S.others}

In this note, we have concentrated on questions of embeddability in
a small number of objects: $\Sym,$ $\Se,$ $\En,$ $\Eq$ and $\Rel,$
with brief observations on a few more: $\Sym/\Sym[_<(\Omega)],$
$(\Rel,\subseteq),$ $(\Rel,\emptyset),$
$\Rel[^{\mathrm{semi}}(\Omega)],$ and objects
which we could write $\Eq[^0(\Omega)],$ $\Eq[^1(\Omega)]$ and
$\Eq[^{0,1}(\Omega)],$ i.e., $\Eq$ regarded as a member of the
variety of lattices with least and/or greatest element.
Similar questions for other objects of the same flavor, for
example the groups of automorphisms of various structures considered in
\cite{B+R}, \cite{MD+RG} and \cite{MD+WCH}, the lattices of congruences
of these objects, etc., would also be of interest.

Above, I quoted results from McKenzie~\cite{mckenzie} only
in the forms in which they were relevant to the questions considered
here; but that paper in fact considers the group $\Sym[(\Omega,\beta)]$
of all permutations of $\Omega$ moving $<\beta$ elements, for a fixed
cardinal $\beta,$ and many of the restrictions proved there
are in terms of $\beta,$ rather than $\cd.$
This, too, represents a direction in which the present results
might be generalized.

The variant of the technique of \S\ref{S.cP+F} illustrated in
Theorem~\ref{T.F(R)}, based on considering functors on
the category $\fb{T.ord}$ of totally ordered sets, rather than on
$\fb{Set},$ also admits wide generalization.
Note that any algebra-valued functor $F$ on $\fb{T.ord}$ can be extended
to the category $\fb{Poset}$ of {\em partially} ordered sets, though
usually not uniquely.
(For instance, by taking each partially ordered
set $P$ to the colimit of the algebras $F(C)$ as $C$ ranges over
all chains in $P;$ or to the limit of the algebras $F(C)$ as
$C$ ranges over all quotient-sets of $P$ given with total orderings
that make the quotient-map isotone.)
Hence in stating Theorem~\ref{T.F(R)} and seeking generalizations,
we could just as well let $F$ be a functor on $\fb{Poset}.$
Still more generally, why not let it be defined on arbitrary
preordered sets?
Or on sets with an arbitrary binary relation?
Or several binary relations?
And for such generalizations, what would be the ``best''
analog of~(\ref{x.genFfin})?
Not knowing what the useful generalizations would be, I have merely
given a sample result.

The referee notes that de~Bruijn's groups with
presentation~(\ref{x.235:2})-(\ref{x.235:5})
are the values of a functor on the category whose objects are
sets $I,$ and whose morphisms are {\em one-to-one} set-maps;
this also applies to several other examples for which we proved
nonembeddability results.
Hence that domain category is {\em not} good from the above
point of view.
A repeated use that we made of non-one-to-one maps in
the arguments of \S\ref{S.cP+F} was in getting left inverses to
one-to-one maps with nonempty domains.
The fact that such inverses exist has the consequence
(also pointed out by the referee) that our key
condition~(\ref{x.genFfin}) is equivalent to the condition
that $F$ respect direct limits.
\vspace{6pt}

A question I have not thought much about, but to which some of
the techniques we have introduced above should be applicable, is

\begin{question}\label{Q.resfin}
What can be said about
groups $A$ embeddable in $(\prod_{n\in\omega}\Sym[(n)])^\Omega,$
monoids $A$ embeddable in $(\prod_{n\in\omega}\Se[(n)])^\Omega,$
$\!K\!$-algebras $A$ embeddable in $(\prod_{n\in\omega}
\En[(K^n)])^\Omega,$ and lattices $A$ embeddable
in $(\prod_{n\in\omega}\Eq[(n)])^\Omega$?

Clearly, such an $A$ must be residually
finite\textup{(}-dimensional\textup{)},
and embeddable in $\Sym,$ $\Se,$ $\En,$ respectively $\Eq.$
Are these conditions on an algebra $A$ sufficient for it to be
embeddable in the indicated object?
\end{question}

I record next a question which it would have been
natural to give in \cite{Sym_Omega:1}, \cite{P_vs_cP}, or
\cite{Sym_Omega:2}, but which occurred to me too late to include in
those papers.
(It is an instance of the direction for investigation suggested in
the second sentence of \cite[\S10]{Sym_Omega:2}.)
There are obvious variants and strengthenings, but for concreteness, I
pose the question here for the object that has been most studied.

\begin{question}\label{Q.Sym=v}
Suppose $G$ and $H$ are subgroups of $\Sym[(\omega)],$ which
together generate that group.
Must $\Sym[(\omega)]$ be finitely generated over one
of these subgroups?
\end{question}

\vspace{6pt}
We end this note with a few items that we have postponed.

\section{Appendix: $\En\cP\En$ can be embedded in
$\En.$}\label{S.jiggleEn}

We shall now give the postponed proof of the above embeddability
statement, Lemma~\ref{L.double}(iii).

Suppose $S$ and $T$ are two nonzero $\!K\!$-algebras (as always,
associative and unital),
and we form their coproduct $S\cP T$ in the variety of such algebras,
calling the coprojection maps $\alpha: S\rightarrow S\cP T$ and
$\beta: T\rightarrow S\cP T.$
Recall \cite[Corollary~8.2]{<>}
that if $B_S,$ $B_T$ are $\!K\!$-vector-space
bases for $S$ and $T,$
containing $1_S$ and $1_T$ respectively, then a $\!K\!$-vector-space
basis for $S\cP T$ is given by the set of finite products
\begin{xlist}\item\label{x.cPbasis}
$\dots\ \alpha(b_i)\ \beta(b_{i-1})\ \alpha(b_{i-2})\ %
\beta(b_{i-3})\ \dots$
\end{xlist}
where, as in~(\ref{x.nmfmg}),
those $b_j$ that are arguments of $\alpha$ (i.e.,
in~(\ref{x.cPbasis}), those with subscript $j$ having the same parity
as $i)$ are taken from $B_S-\{1_S\},$ those that are
arguments of $\beta$ are taken from $B_T-\{1_T\},$ and
the empty product is understood to give the identity element
$1=\alpha(1_S)=\beta(1_T).$
As in~(\ref{x.nmfmg}), I have not shown the first
and last terms, since each may be either
an $\!\alpha\!$-term or a $\!\beta\!$-term.

Letting $S=T=\En,$ and letting $B$ be a $\!K\!$-basis for this
algebra containing~$1,$ a basis for $\En\cP\En$ is thus given by
the words~(\ref{x.cPbasis}) with all $b_j$ taken from $B-\{1\}.$
To show $\En\cP\En$ embeddable in $\En,$ it will suffice to find
a representation of $\En\cP\En$ by $\!K\!$-linear endomorphisms on
a vector space of the same dimension as $V,$ such that the images
of the elements~(\ref{x.cPbasis}) are linearly independent.
For this, in turn, it will suffice to find a family of $\leq\cd$
representations of $\En\cP\En$ on such spaces, such that every
nontrivial linear relation among elements~(\ref{x.cPbasis}) fails to
hold in at least one of these representations, since then all such
relations will fail in the direct sum of the representations.

The representations we use will each be on
$\bigoplus_\omega V,$ a direct sum of countably
many copies of $V,$ with basis $\Omega\times\omega.$
For each $k\in\Omega$ we shall call the $\!k\!$th copy of $V,$
i.e., the span of $\Omega\times\{k\},$ the ``$\!k\!$th level''
of $\bigoplus_\omega V,$ and, extending the notation we are using
on its basis, we shall denote the element at the $\!k\!$th
level corresponding to any $v\in V$ by $(v,k).$
We define the {\em natural action} of each $f\in\En$ on
the $\!k\!$th level of $\bigoplus_\omega V$ to be given by
$f(v,k)=(f(v),k),$ i.e., to mimic its action on $V.$
The natural action on $\bigoplus_\omega V$
will mean the direct sum of these actions.

We now define our $\cd$ actions of $\En\cP\En$ on $\bigoplus_\omega V.$
They will be indexed by the set of those vector space automorphisms $t$
of $\bigoplus_\omega V$ which have order $2,$ fix all but
finitely many members of the basis $\Omega\times\omega,$ and take
the remaining members of that basis to linear combinations
of members of that basis with integer coefficients.
The last two conditions insure there are only $\cd$ such $t.$
For each such choice of $t,$ we map $\En\cP\En$
to $\En[(\bigoplus_\omega V)]$ by sending each element
$\alpha(f)$ $(f\in\En)$ to the natural action of $f$
on $\bigoplus_\omega V,$ which we denote by the
same symbol $f,$ while we send $\beta(f)$ to $t\,f\,t^{-1}.$

For the remainder of the proof, let us fix a nonzero $x\in\En\cP\En.$
If $x\in K,$ then clearly $x$ has nonzero action under all
of our representations; so let us assume $x\notin K,$ and show how
to construct a $t$ such that the action of $x$ on $\bigoplus_\omega V$
under the representation indexed by $t$ is nonzero.

To do this, let us fix an arbitrary element $r\in\Omega,$ and note
that $\En$ will be the direct sum of the \mbox{$\!1\!$-dimensional}
subspace spanned by the identity map, and the subspace $\En_0$ of
consisting of those maps that take $r$ to a linear combination
of elements of $\Omega-\{r\}.$
Choosing, temporarily, an arbitrary basis $B'$ of $\En$ consisting
of $1$ and elements of $\En_0,$ let us express $x$ as a linear
combination of words of the form~(\ref{x.cPbasis}) with the $b_j$
taken from $B'-\{1\}.$
Let $n$ be the maximum of the lengths
of the words occurring with nonzero coefficient in this expression.
Let $S\subseteq B'$ be the set consisting of $1$ and all
those elements of $B'-\{1\}$ that occur
(as arguments of $\alpha$ or $\beta)$ in the expression for $x.$
From the fact that the span of $S$ is a finite dimensional
subspace of $\En,$ it is not hard to see that there will exist a
finite subset $\Sigma\subseteq\Omega$ containing the
element $r,$ such that, if we write
$P_\Sigma\in\En$ for the element that fixes all members
of $\Sigma$ and annihilates all members of $\Omega-\Sigma,$ then
the linear operator on $\En$ given by
\begin{xlist}\item\label{x.PfP}
$f\ \mapsto\ P_\Sigma\,f\,P_\Sigma$
\end{xlist}
is one-to-one on the span of $S.$
Let us fix such a set $\Sigma\subseteq\Omega\<.$

We can now describe the basis $B$ of $\En$ in terms of which
we will work for the rest of the proof.
For all $p,q\in\Omega,$ let $E(q,p)\in\En$ be the linear map that
takes $p$ to $q,$ and all other members of $\Omega$ to $0.$
Since the typical member of $\En$ has infinite-dimensional
range, the elements $E(q,p)$ do not span $\En;$ but
their linear combinations do give all possible behaviors on our
finite set $\Sigma,$ which is what we will need.
Let us choose $B$ to consist of the identity operator $1,$
all the operators $E(q,p)$ with $p\in\Sigma$
(and $q$ unrestricted) except for $E(r,r),$
and the members of any basis of the space
of those endomorphisms that annihilate $\Sigma.$

To see that an arbitrary $f\in\En$ may be represented by a linear
expression in members of $B,$ first set the coefficient of~$1$ in
this expression to be the coefficient of $r$ in $f(r).$
Subtracting from $f$ that multiple of $1$ gives a member of $\En_0,$
whose behavior on the elements of $\Sigma$ can be represented by
a finite linear combination of the operators $E(q,p)$ with
$p\in\Sigma$ and $(p,q)\neq(r,r).$
Subtracting this off, we are left with an operator annihilating
$\Sigma,$ which can be uniquely represented using the
elements introduced in the last part of our definition of $B.$
Clearly, this expression for $f$ is unique.

We now take our earlier expression for $x$ as a linear
combination of words~(\ref{x.cPbasis}) with all $b_j\in B'-\{1\},$
and substitute for the $b_j$ their
expressions as linear combinations of elements of $B-\{1\},$ getting
an expression for $x,$ again as a combination of
words~(\ref{x.cPbasis}), with the $b_j$ now in $B-\{1\}.$
Clearly, these words still have length $\leq n.$
I claim, moreover, that by our choice of $\Sigma,$ the expression
contains, with nonzero coefficient, at least one
product~(\ref{x.cPbasis}) of length $n$ in which all $b_j$
have the form $E(q,p)$ with $p$ and $q$ both in $\Sigma$ (and
by definition of $B,$ with $(p,q)\neq(r,r)).$

To see this, note that a consequence of our normal form
for coproducts of associative $\!K\!$-algebras is that
$\En\cP\En$ can be identified with a direct sum of iterated
tensor products $\En_0\,\otimes_K\dots\otimes_K\,\En_0,$ where the
$\!0\!$-fold tensor product, i.e., $K,$ occurs once,
and each higher tensor product occurs twice,
corresponding to the two ways of labeling the tensor factors
alternately with $\alpha$ and $\beta.$
This identification maps any element of one of these direct summands
(e.g., an element $v_1\otimes v_2\otimes v_3$ in the summand
labeled with $\alpha,\ \beta,\ \alpha)$ to the corresponding
sum of products of elements of $\alpha(\En_0)$ and $\beta(\En_0)$
(e.g., $\alpha(v_1)\<\beta(v_2)\<\alpha(v_3)$ in that example),
and is therefore independent of choice of basis.
(It does depend on our choice of linear complement $\En_0$
for $K$ in $\En,$ which depends on our choice of $r,$
but we made that choice at the start and have not changed it.)
I claim that the condition by which we chose the finite set
$\Sigma$ implies that if the $\!K\!$-linear map~(\ref{x.PfP})
is applied simultaneously to every tensor factor $\En_0$ in every
summand in the above expression for $\En\cP\En,$ the element $x$
continues to have nonzero components in all the
degrees where $x$ had them; in particular, in degree $n.$
Indeed, the linear relations holding among a set of expressions in a
tensor product of vector spaces depend only on the linear relations
holding among the elements of the given spaces that occur
in these expressions; and $\Sigma$ was chosen
so that~(\ref{x.PfP}) creates no new linear relations among
the elements occurring in our original expression for $x.$
Now in terms of our new basis $B-\{1\}$ of $\En_0,$ the
map~(\ref{x.PfP}) acts by throwing out all basis elements other
than the $E(q,p)$ with $p,q\in\Sigma.$
So since $x$ continues to have nonzero component in degree $n$
after the application of~(\ref{x.PfP}),
the expression for $x$ using the basis $B-\{1\}$ does indeed involve
at least one length-$\!n\!$ word in such elements $E(q,p)$ alone.

We now choose, subject to a restriction to be given in a moment,
a particular length-$\!n\!$ word of this sort
occurring with nonzero coefficient in $x,$
\begin{xlist}\item\label{x.Eprod}
$\dots\ \alpha(E(q_i,p_i))\ \beta(E(q_{i-1},p_{i-1}))\ %
\alpha(E(q_{i-2},p_{i-2}))\ \beta(E(q_{i-3},p_{i-3}))\ \dots\,,$\\
where $(q_j,p_j)\in\Sigma\times\Sigma-\{(r,r)\}$ $(j=1,\dots,n).$
\end{xlist}
To do this, let us work from the right, first choosing
an $\alpha(E(q_1,p_1))$ or $\beta(E(q_1,p_1))$
that occurs as a rightmost factor in some length-$\!n\!$ word of the
form~(\ref{x.Eprod}) in our expression for $x;$ then choosing for
$\beta(E(q_2,p_2))$ or $\alpha(E(q_2,p_2)),$ as the case may be,
an element of this form that
occurs in second position from the right, immediately to the left
of our first chosen factor, in at least one
length-$\!n\!$ word of that form; and so on.
The one restriction we impose is
that for each $j=1,\dots,n,$ $(q_j,p_j)$ should satisfy
$q_j\neq p_j$ if this is possible, i.e., if there
{\em is} a factor satisfying $q_j\neq p_j$ which occurs,
followed by the terms chosen so far, in the $\!j\!$th position of a
length-$\!n\!$ word occurring in $x.$
Henceforth,~(\ref{x.Eprod}) will denote the particular word so chosen.

The reason we avoid choices with $q_k=p_k$ whenever possible
is that it will not be as easy to make use of the
fact that elements of the form $E(p,p)$ are nonscalar as it
will for other elements; but we will be able to do so if there are no
elements $E(q,p)$ with $q\neq p$ ``in the vicinity''.
The distinction between these cases is used in the next definition.

For $k=1,\dots,n,$ let us define elements $p_k'\in V,$ by
\begin{xlist}\item\label{x.p'}
$\!p_k'=p_k\in\Omega$\quad\quad if $q_k\neq p_k,$\\[2pt]
$p_k'=p_k+r\in V$ if $q_k=p_k.$
\end{xlist}
Note that in the second line above, the summands $p_k$ and $r$
are distinct members of $\Omega,$ for
if $q_k=p_k,$ we cannot have $p_k=r,$ since $E(r,r)\notin B.$
For those values of $k$ such that $q_k=p_k,$ let us temporarily
form a new basis of the $\!k\!$th level (the $\!k\!$th direct
summand) of $\bigoplus_\omega V,$ by deleting
from $\Omega\times\{k\}$ the basis element $(r,k)$
and inserting $(p_k',k)$ in its place, while for those $k$ such
that $q_k\neq p_k,$ let us keep the basis $\Omega\times\{k\}.$
Thus, the union over $k$
of these bases is a new basis of $\bigoplus_\omega V.$
Note that for {\em every} $k,$
$(p_k',k)$ and $(q_k,k)$ are distinct elements of our new basis.

We now define our automorphism $t$ of $\bigoplus_\omega V$
to fix all elements of our new basis, except for the following
$2(n+1)$ elements, which we let it transpose in pairs as shown:
\begin{xlist}\item\label{x.Etranspose}
$(p_1,0)\leftrightarrow (p_1',1),$\quad
$(q_{k-1},k{-}1)\leftrightarrow (p_k',k)$ $(1<k\leq n),$\quad
$(q_n,n)\leftrightarrow (q_n,n{+}1).$
\end{xlist}

This is where we need the second line of~(\ref{x.p'}).
It insures that even if $q_k=p_k,$ the basis elements $(p_k',k)$
and $(q_k,k)$ are distinct, so
that the rules $(q_{k-1},k{-}1)\leftrightarrow (p_k',k)$ and
$(q_k,k)\leftrightarrow (p_{k+1}',k{+}1)$ do
not give contradictory specifications of the action of $t$
on the same basis element at the $\!k\!$th level.

As promised, $t$ has order $2$ and sends members of our original basis
$\Omega\times\omega$ to linear combinations of members of that basis
with integer coefficients.
Having noted this, we shall now go back to using
$\Omega\times\omega$ as the basis in terms of which we will
compute with elements of $\bigoplus_\omega V,$ and shall
think of the terms of~(\ref{x.Etranspose}) as
expressions in those basis elements.

As sketched earlier, we now
let $h: \En\cP\En\rightarrow\En[(\bigoplus_\omega V)]$
be the homomorphism
which takes elements $\alpha(f)$ to $f$ (acting by the natural
action) and elements $\beta(f)$ to $t\,f\,t^{-1}=t\,f\,t.$

The remainder of the proof follows closely the concluding
steps of the proof of Lemma~\ref{L.double}(ii).
We want to show that $h(x)\neq 0.$
Clearly, this is equivalent to proving nonzero the element
$h(x)'$ that we get on multiplying $h(x)$ on the right by $t$
if the rightmost term of~(\ref{x.Eprod}) is $\alpha(E(p_1,q_1)),$
while leaving that side unchanged if that term
is $\beta(E(p_1,q_1)),$ and multiplying on the left by $t$
if the leftmost term of~(\ref{x.Eprod}) is $\alpha(E(p_n,q_n)),$
while leaving that side unchanged if that term is $\beta(E(p_n,q_n)).$
To avoid cumbersome language, we shall call
products of $\!E(p,q)\!$'s and $\!t\!$s occurring with
nonzero $\!K\!$-coefficient in our expression for $h(x)'$ the
``summands'' in that expression (not counting the $\!K\!$-coefficients
as parts of these ``summands'').
A consequence of our definition of $h(x)'$ is that the summand therein
arising from the term~(\ref{x.Eprod}) of $h(x)$ has a $t$ at each
end, so we can now write it so as show those ends.
Let us also give it a name:
\begin{xlist}\item\label{x.tEprod}
$u\ =\ t\ E(q_n,p_n)\ t\ E(q_{n-1},p_{n-1})\ t\ \dots\ %
t\ E(q_2,p_2)\ t\ E(q_1,p_1)\ t.$
\end{xlist}
The general summand in the expression for $h(x)'$ will likewise
be an alternating product of $t$ and elements
of $B-\{1\}$ (which we will call the ``$\!B\!$-factors'' of the
summand), with at most $n$ of the latter.
We claim now that when we apply $h(x)'$ to $(p_1,0),$ the summand
$u$ shown in~(\ref{x.tEprod}), and only that summand, leads to a
nonzero component at the $\!n{+}1\!$st level of $\bigoplus_\omega V.$

Since in each summand, the only factors that
carry elements of $\bigoplus_\omega V$ from one level to the next are
the factors $t,$ and each of these moves elements by only one level,
the only summands in $h(x)'$ that can possibly lead to components
at the $\!n{+}1\!$st level in $h(x)'(p_1,0)$ are those which,
like~(\ref{x.tEprod}), have exactly $n$ (rather than fewer)
$\!B\!$-factors, and have a $t$ at each end.
Consider any such summand
\begin{xlist}\item\label{x.w}
$w\ =\ t\ b_n\ t\ b_{n-1}\ t\ \dots\ t\ b_2\ t\ b_1\ t,$
\end{xlist}
where $b_k\in B-\{1\}$ $(k=1,\dots,n).$
We shall show inductively for $k=1,\dots,n{+1}$
that if we apply to $(p_1,0)$ the substring
$t\,b_{k-1}\,t\,b_{k-2}\,t\,\dots\,t\,b_2\,t\,b_1\,t$ of $w,$ then the
components of the result in levels higher than the $\!k\!$th
are zero, the component in the $\!k\!$th level is a scalar multiple
of $(p_k',k),$ and the scalar factor is nonzero if and
only if our substring agrees exactly with the corresponding
substring of $u,$ i.e., if and only if $b_i=E(q_i,p_i)$ for
$i=1,\dots,k-1.$

The base case $k=1$ is immediate: we are merely applying $t$
to $(p_1,0),$ and by~(\ref{x.Etranspose}) the result is $(p_1',1).$
(This rightmost factor $t$ was important in distinguishing
the action of $u$ from actions of summands of $h(x)'$ not ending
in $t;$ but we have already used it to exclude such strings
from consideration.)

Now let our inductive assumption hold for some $k$ with $1\leq k\leq n.$
If we do not have $b_i=E(q_i,p_i)$ for all $i<k,$ then by that
inductive assumption, $t\,b_{k-1}\,t\,\dots\,t\,b_1\,t\,(p_1,0)$
has zero component at the $\!k\!$th and higher levels, so
multiplication by $b_k$ and then by $t$ will not bring anything to
the $\!k{+}1\!$st level or higher.

If the factors of $w$ so far {\em have}
agreed with those of $u,$ then by inductive hypothesis,
$t\,b_{k-1}\,t\,\dots\linebreak[0]\,t\,b_1\,t\,(p_1,0)$
has at the $\!k\!$th level a nonzero scalar multiple of $(p_k',k).$
When we apply $b_k$ to this, if $b_k=E(q_k,p_k)$
i.e., if this too agrees with the corresponding factor
of $u,$ then we see that, whether $(p_k',k)$ is equal
to $(p_k,k)$ or to $(p_k,k)+(r,k),$ the factor $E(q_k,p_k)$ will send
it to $(q_k,k);$ and the subsequent application of $t$
will give us a term $(p_{k+1}',k{+}1)$ at the $\!k{+}1\!$st level,
as desired.

If $b_k\neq E(q_k,p_k),$ there are several cases to consider.
First, $b_k$ might be one of the members of $B$ belonging to
the subspace of $\En_0$ annihilating $\Sigma.$
In that case, it annihilates $(p_k',k),$
leaving nothing at the $\!k\!$th level, and the subsequent
application of $t$ brings nothing to the $\!k{+}1\!$st level.
Otherwise, we have $b_k=E(q,p)$ for some $(q,p)\neq(q_k,p_k).$
Clearly, the only cases in which $E(q,p)$ can fail to annihilate
$(p_k',k)$ are

(a)~if $p=p_k,$ or

(b)~if $p_k'$ has the form $p_k+r,$ and $p=r.$

In case~(a), since we have assumed that $E(q,p)\neq E(q_k,p_k),$
we must have $q\neq q_k,$ so $E(q,p)(p_k',k)=(q,k)\neq(q_k,k),$
and by~(\ref{x.Etranspose}), the subsequent application of $t$ will
not bring this up to the $\!k{+}1\!$st level.

In case~(b), the assumption $p_k'=p_k+r$ means, by~(\ref{x.p'}),
that $E(q_k,p_k)=E(p_k,p_k).$
But recall that in choosing the term~(\ref{x.Eprod}), we avoided
this possibility whenever possible.
A consequence is that since $w$ was not chosen in preference to $u,$
we must likewise have $b_k=E(p,p).$
However, case~(b) assumed $p=r,$ so $b_k=E(r,r),$
which is excluded by the definition of $B.$
Hence case~(b) does not occur.

We thus conclude that $u$ is the unique summand in $h(x)'$ which,
when applied to $(p_1,0),$ gives an element having nonzero
component at the $\!n{+}1\!$st level.
Hence $h(x)'\,(p_1,0)\neq 0,$ so $h(x)'\neq 0,$ so $h(x)\neq 0,$
completing the proof that our action of $\En\cP\En$ is
faithful.\qed\vspace{6pt}

Now for an unexpected bonus.
We saw in \S\ref{S.notEn} that when $K$ had large cardinality, our
tools for proving {\em nonembeddability} of algebras in $\En$ gave
weaker results than we expected from our results on groups,
monoids and lattices.
We shall now see that we also get stronger {\em positive}
embeddability results in that situation.
We need the following fact.

\begin{lemma}[cf.\ {\cite[Exercise~1, p.248]{NJ}}]\label{L.kG}
Let $k$ be any commutative integral domain, and $M$ the multiplicative
monoid of $k.$
Then the product $\!k\!$-algebra $k^\omega$ contains a subalgebra
isomorphic to the monoid algebra $k\<M.$
In particular, letting $G$ denote the group of units of $k,$
it contains a copy of the group algebra~$k\<G.$
\end{lemma}\begin{proof}
For each $a\in k,$ let $x_a\in k^\omega$ denote the
sequence of powers, $(1,a,a^2,\dots).$
Clearly, $a\mapsto x_a$ is a monoid homomorphism
from $M$ into the multiplicative monoid of $k^\omega.$
Moreover, the elements $x_a$ are $\!k\!$-linearly independent
by the properties of the Vandermonde determinant.
Hence they span a subalgebra of $k^\omega$ isomorphic to $k\<M.$
The final assertion clearly follows.
\end{proof}

\begin{proposition}\label{P.cP_K}
$\En$ contains a coproduct of $\cd[K]$ copies
of itself as an associative $\!K\!$-algebra.
\textup{(}Hence, in view of Theorem~\ref{T.bigcP}\textup{(iii)},
it contains such a coproduct of
$\max(2^{\cd},\linebreak[0]\cd[K])$ copies of itself.\textup{)}
\end{proposition}\begin{proof}
For $K$ finite, our main statement is weaker
than Theorem~\ref{T.bigcP}(iii), so assume $K$ infinite.

Since $\En$ contains a subalgebra isomorphic to $K^\omega,$
the preceding lemma shows that it contains a copy of the group
algebra $K\<G$ on a group $G$ of cardinality $\cd[K-\{0\}]=\cd[K].$

We have also just seen that it contains a copy of $\En\cP\En;$
combining these observations, we conclude
that it contains a copy of $\En\cP K\<G.$
For notational convenience, let us identify $\En$ with the first
factor in this coproduct, and $G$ (which we regard as an abstract
group, forgetting about its relation with $K)$ with its image
in the second factor.
Then for each $g\in G,$ we can form the conjugate algebra
$g~\En~g^{-1},$ getting $\cd[K]$ isomorphic copies of $\En.$
Moreover, from the linear independence of the elements
of $G$ in $K\<G$ and the normal
form for the coproduct $\En\cP K\<G,$ one sees that
the subalgebra these generate will be their coproduct, so $\En$
indeed contains a coproduct of $\cd[K]$ copies of itself.
(This is the ring-theoretic analog of the appearance of big free
products inside smaller free products in the Kurosh Subgroup Theorem.)
The final assertion clearly follows.
\end{proof}

\begin{corollary}\label{C.cP_K}
$\En$ contains a coproduct, as $\!K\!$-algebras, of
$\max(2^{\cd},\cd[K])$ copies of every $\!K\!$-algebra of dimension
$\leq\cd.$\qed
\end{corollary}

Of course, these results lead to the questions:

\begin{question}\label{Q.En_big}
Can we prove a version of Theorem~\ref{T.funct'}\textup{(iii)}
with $F(2^{\cd})$ replaced by
$F(\max(2^{\cd},\linebreak[0]\cd[K]))$ for a wider class of
functors $F$ than the ``coproduct as $\!K\!$-algebras of $I$ copies
of $A$'' functors of Corollary~\ref{C.cP_K}?
\end{question}

An interesting test case would be that in which $F$ is the functor
taking $I$ to the $\!K\!$-algebra presented
by an $\!I\!$-tuple of commuting idempotents.

Remark: One can strengthen Proposition~\ref{P.cP_K} so as to increase
$\max(2^{\cd},\linebreak[0]\cd[K])$ to
$\cd[K]^{\cd},$ which occasionally exceeds the former value; e.g.,
when $\cd[K]$ has the form $\aleph_{\alpha+\omega}$ and is $>2^{\cd}.$
To do so, replace $\omega$ in Lemma~\ref{L.kG} with the free
abelian monoid $A(\Omega)$ on $\Omega,$ which has cardinality $\cd,$
replace the maps $x_a$ in the proof by all homomorphisms of
$A(\Omega)$ into the multiplicative monoid $M$ of $k,$ apply the
theorem on linear independence of characters to conclude that
$k^{A(\Omega)}$ contains a copy of the monoid algebra $k\,M^\Omega,$
and use this version of the lemma to get the strengthened
proposition.
(One can't raise the cardinal in the proposition higher
than $\cd[K]^{\cd},$ since $\dim_K(\En)=\dim_K(V^\Omega)\leq
\dim_K((K^\Omega)^\Omega)=\dim_K(K^{\Omega\times\Omega})=
\dim_K(K^\Omega)\leq\cd[K^\Omega]=\cd[K]^{\cd};$ cf.\ \cite[p.247,
Theorem~2]{NJ}.)

It is curious that the nonembeddability results we are
able to prove become weaker than expected as soon as $\cd[K]>\cd,$ but
our positive results become stronger only when $\cd[K]>2^{\cd}.$

\section{Appendix: another embedding of $\Sym\cP\Sym$
in $\Sym.$}\label{S.paths}

Here is the alternative proof of Lemma~\ref{L.double}(i) mentioned
shortly before the statement of Lemma~\ref{L.double}.

As in the preceding appendix, we begin by recalling a
structure theorem for coproducts, this time coproducts of groups.
But we will make different assumptions (we won't assume there are
only two groups, but we will assume the groups are
disjoint except for their identity elements, so that we do
not have to write the coprojection maps explicitly), we will
use the structure theorem, initially, for a different purpose
(to motivate a somewhat bizarre action of the coproduct group), and we
will consider for most of this section arbitrary groups (or sometimes,
monoids) rather than a pair of copies of $\Sym,$ though that is the
case to which we will ultimately apply our result.

Let $(G_i)_{i\in I}$ be any family of groups whose sets of
nonidentity elements are pairwise disjoint.
Recall that the general element of
the coproduct group $\coprod_{i\in I} G_i$ can be written uniquely as
\begin{xlist}\item\label{x.gng1}
$g_n\ \dots\ g_2\ g_1,$
\end{xlist}
where $n\geq 0$ and $g_n,\dots, g_1\in\bigcup_{i\in I} G_i-\{1\},$
say with
\begin{xlist}\item\label{x.i_r}
$g_r\in G_{i_r}$ $(r=1,\dots,n),$
\end{xlist}
and where successive
indices $i_r,\ i_{r+1}$ are distinct; i.e., two elements from the
same group $G_i$ never occur in immediate succession.

Suppose we are given a representation of each $G_i$ by permutations of a
nonempty set $\Omega_i,$ and we form the direct product $\prod\Omega_i.$
Then starting at any point $(x_i)_{i\in I}\in\prod\Omega_i,$ an
expression~(\ref{x.gng1}) allows us to construct
a ``path'' in $\prod\Omega_i:$
At the first step, we move from $(x_i)$ to the point agreeing
with $(x_i)$ on all but the $\!i_1\!$-coordinate, which
is moved by $g_1;$ the next step takes us to
the point whose $\!i_2\!$-coordinate has also been moved by $g_2,$ etc..
This suggests that we make the set of ``paths'' in $\prod\Omega_i$ in
which each step involves changing just one coordinate, and two
successive
steps never change the same coordinate, into a $\!\coprod G_i\!$-set.

Some difficulties arise.
Though, by assumption, no $g_r$ in~(\ref{x.gng1}) is an
identity element, some of these factors may lie in the
stabilizers of the coordinates they are to be applied to.
In such cases should we, as the above description might suggest,
allow trivial ``steps'' in our path, where no coordinate is changed?
It turns out that this would lead to difficulties; so we
specify that in such cases, no step is added to our path.
Also, for inverses to behave correctly, we must allow some elements
to delete rather than adding steps to our paths.

The resulting construction is described in the next lemma.
Note that it is not claimed that if the $\Omega_i$ are faithful
$\!G_i\!$-sets, then the
$\!\coprod G_i\!$-set $\boxtimes_{i\in I}\Omega_i$
described there is also faithful.
Rather, we shall subsequently note additional conditions that ensure
faithfulness.

Up to the step of achieving faithfulness, our
construction works as well for arbitrary
monoids as for groups, so the lemma below is stated in that context.
As usual, we understand actions to be {\em left} actions.
Note that in the statement below, subscripts $n$ and $r$ do not
correspond to the subscript $i$ above; rather, each of the steps
$x_r$ $(r=1,\dots,n)$ in our ``path'' is itself
an $\!I\!$-tuple $x_r=(x_{r,i})_{i\in I}.$

\begin{lemma}\label{L.boxX}
Let $(M_i)_{i\in I}$ be a family of monoids, and for each $i\in I$
let $\Omega_i$ be an $\!M_i\!$-set.
Let $\boxtimes_{i\in I} \Omega_i$ be the set of all finite sequences
$(x_1,\dots,x_n)$ $(n\geq 1)$ where each $x_r\in\prod_{i\in I} \Omega_i$
$(r=1,\dots,n),$
every pair of successive terms $x_r,\ x_{r+1}$ $(1\leq r<n)$
differs in one and only one coordinate, and the coordinate at which
$x_{r+1}$ differs from $x_r$ is not the same as the coordinate
at which $x_r$ differs from $x_{r-1}$ $(1<r<n).$

Then $\boxtimes_{i\in I} \Omega_i$ may be made a
$\coprod_{i\in I} M_i\!$-set \textup{(}which we might call the
``path product'' of the $\!M_i\!$-sets $\Omega_i)$
by defining the action of each $g\in M_j$ $(j\in I)$ on $x=
(x_1,\dots,x_n)\in\boxtimes_{i\in I} \Omega_i$ by the following rules:
\begin{xlist}\item\label{x.action}
Let $x_n'\in\prod_{i\in I} \Omega_i$ denote the $\!I\!$-tuple obtained
from $x_n$ by modifying its $\!j\!$-coordinate via the action of~$g$
on that element of $M_j.$
Then\\
\textup{(i)}~ If $x_n'=x_n,$ we define $gx=x.$\\
\textup{(ii)}~
If $x_n'\neq x_n,$ then:\\
\textup{(ii.a)}~
If either $n=1,$ or if $n>1$ and the coordinate at which
$x_n$ differs from $x_{n-1}$ is not the $j\!$th, we define
$gx=(x_1,\dots,x_n,x_n'),$\\
\textup{(ii.b)}~ If $n>1$ and the coordinate at which
$x_n$ differs from $x_{n-1}$ {\em is} the $j\!$th, then:\\
\textup{(ii.b.1)}~ If $x_n'\neq x_{n-1},$
we define $gx=(x_1,\dots,x_{n-1},x_n'),$ while\\
\textup{(ii.b.2)}~ If $x_n'=x_{n-1},$ we define
$gx=(x_1,\dots,x_{n-1}).$
\end{xlist}
\end{lemma}\begin{proof}
By the universal property of the coproduct $\coprod_{i\in I} M_i,$ an
action of that monoid on the set $\boxtimes_{i\in I} \Omega_i$ will
be defined if we verify that for each $j\in I,$ the above conditions
define an action of $M_j$ on that set.

This could be done by brute force, dividing into cases according
to which headings of the above definition the actions of {\em two}
successive elements of $M_j$ come under.
But there is a trick that greatly simplifies this calculation
(cf.\ \cite[proof of Proposition~3.6.5]{245}).
For each $j,$ we shall define a bijection $\varphi_j$ between
$\boxtimes_{i\in I} \Omega_i$ and a set
$(\boxtimes_{i\in I} \Omega_i)_{(j)},$ such
that it will be easy to define an action of
$M_j$ on $(\boxtimes_{i\in I} \Omega_i)_{(j)},$ and also easy to verify
that when we transport this action from
$(\boxtimes_{i\in I} \Omega_i)_{(j)}$ to $\boxtimes_{i\in I} \Omega_i$
via the bijection $\varphi_j,$ the resulting action
is described by~(\ref{x.action}).
Roughly, elements of $(\boxtimes_{i\in I} \Omega_i)_{(j)},$ like
elements of $\boxtimes_{i\in I} \Omega_i,$ will represent ``paths'' in
$\prod\Omega_i,$ but in $(\boxtimes_{i\in I} \Omega_i)_{(j)}$ we
require every
such path to have a final step involving the $\!j\!$-coordinate,
at the price of allowing this step (and only this step) to be trivial
(i.e., to satisfy $x_{n-1}=x_n).$

Here is the precise description.
We take the elements of $(\boxtimes_{i\in I} \Omega_i)_{(j)}$ to be
sequences $x=(x_1,\dots,x_n)$ of elements of $\prod\Omega_i,$ this
time with $n$ always $\geq 2,$ such that

(i)~ $(x_1,\dots,x_{n-1})\in\boxtimes_{i\in I} \Omega_i,$

(ii)~ $x_{n-1}$ and $x_n$ are either equal,
or differ in the $\!j\!$th coordinate only.

(iii)~ If $n>2,$ the coordinate
at which $x_{n-2}$ and $x_{n-1}$ differ is {\em not} the $\!j\!$th.

Note that none of these conditions constrains the $\!j\!$th
coordinate of $x_n.$
Hence a monoid action of $M_j$ on $(\boxtimes_{i\in I} \Omega_i)_{(j)}$
may be defined by the rule
\begin{xlist}\item\label{x.gx1xn}
$g\,(x_1,\dots,x_{n-1},x_n)\ =\ (x_1,\dots,x_{n-1},x_n'),$
\end{xlist}
where, as in the statement of the lemma, $x_n'$ is obtained from
$x_n$ by modifying its $\!j\!$th coordinate by the action of $g,$
and leaving all other coordinates unchanged.

Let us now define
$\varphi_j:\boxtimes_{i\in I} \Omega_i\rightarrow(\boxtimes_{i\in I}
\Omega_i)_{(j)}$ to leave unchanged all
$x=(x_1,\dots,x_n)\in\boxtimes_{i\in I} \Omega_i$ with $n>1$ in which
the coordinate at which $x_{n-1}$ and $x_n$ differ
happens to be the $\!j\!$th, while appending to every other element
of $\boxtimes_{i\in I} \Omega_i$ a repetition of its final term.
It is straightforward that this is a bijection, and easy to verify
that the action of $M_j$ on $\boxtimes_{i\in I} \Omega_i$ induced, via
this bijection, by its action~(\ref{x.gx1xn})
on $(\boxtimes_{i\in I} \Omega_i)_{(j)},$ is as described
in~(\ref{x.action}), completing our proof.
\end{proof}

We now want to use the above construction to get {\em faithful} actions.
We shall see that we can do this if our monoids are right cancellative,
in particular, if they are groups.

Let us call an action of a monoid $M$ on a set
$\Sigma$ {\em strongly faithful} if for every finite family
of distinct elements $g_1,\dots,g_n\in M,$ there exists
$y\in \Sigma$ such that $g_1y,\dots,g_ny$ are distinct.
(The $n=0$ case of this condition says that $\Sigma$ is nonempty.)
It is easy to see that if $\Sigma$ is a faithful
$\!M\!$-set, then the disjoint union $1\cup\Sigma\cup\Sigma^2\cup
\Sigma^3\cup\dots$ (indexed by the
natural numbers, with each $\Sigma^n$ given the
product $\!M\!$-set structure) is strongly faithful, and that
this construction does not increase infinite cardinalities.
On the other hand, we have

\begin{lemma}\label{L.fthfl}
Let $(M_i)_{i\in I}$ be a family of monoids, each having the right
cancellation property
\begin{xlist}\item\label{x.hg=g}
$ac=bc\ \implies\ a=b$ $(a,\,b,\,c\in M_i),$
\end{xlist}
and for each $i\in I,$
let $\Omega_i$ be a {\em strongly} faithful left $\!M_i\!$-set.
Then $\boxtimes_{i\in I} \Omega_i,$ defined and made a
$\coprod_{i\in I} M_i\!$-set
as in Lemma~\ref{L.boxX}, is a faithful $\coprod_{i\in I} M_i\!$-set.
\end{lemma}\begin{proof}
Given distinct elements $g=g_m\dots g_2\<g_1$ and
$h=h_n\dots h_2\<h_1$ in $\coprod M_i,$ we wish to
find an $x\in\boxtimes_{i\in I} \Omega_i$ such that $gx\neq hx.$
We shall construct below such an element which is of length 1, i.e., a
$\!1\!$-tuple $(x_1)$ with $x_1\in\prod \Omega_i.$
To do this let us, for each $j\in I,$ choose the
$\!j\!$-coordinate of $x_1$ as follows.
Let $r_{j,1}<\dots<r_{j,m_j}$ be the values of $r$ for
which $g_r\in G_j$
and $s_{j,1}<\dots<s_{j,n_j}$ the values $s$ for which $h_s\in G_j.$
Although the elements
\begin{xlist}\item\label{x.gr1grm}
$1,\ g_{r_{j,1}},\ g_{r_{j,2}}g_{r_{j,1}},\ \dots,\ %
g_{r_{j,m_j}}...\,\<g_{r_{j,2}}g_{r_{j,1}},$
\quad and \quad
$1,\ h_{s_{j,1}},\ h_{s_{j,2}}h_{s_{j,1}},\ \dots,\ %
h_{s_{j,n_j}}...\,\<h_{s_{j,2}}h_{s_{j,1}}$
\end{xlist}
of $M_j$ need not all be distinct, {\em successive} elements of
each of these lists will be so, by right cancellation.
For each $j\in I$ let us use strong faithfulness of $M_j$ to choose, as
the $\!j\!$-coordinate of $x_1,$ an element
of $\Omega_j$ whose images under {\em distinct} elements of
the combined list~(\ref{x.gr1grm}) are distinct.
(For some $j,$ the $m_j$ and $n_j$ in~(\ref{x.gr1grm}) may
both be zero, in which case this condition is vacuous.)
Thus, the images of $x_1$ under successive terms of
each list in~(\ref{x.gr1grm}) are distinct.

It follows from this choice that each time we apply an element $g_r$
or $h_r$ in building up $g_m\dots g_1 x$ or
$h_n\dots h_1 x,$ it in fact moves the coordinate
of the term of $g_{r-1}\dots g_1 x,$ respectively
$h_{r-1}\dots h_1 x$ to which it is applied, and successive steps move
different coordinates.
Thus, we are always in case~(ii.a) of the
definition of the action of $M_{i_r},$
so the final elements $g_m\dots g_1 x$ and $h_n\dots h_1 x$
reflect all the steps of this process.
In particular, if $m\neq n,$ the elements $gx$ and $hx$ have
different lengths, if $m=n$ and the sequence of indices
in $I$ determined by $g$ and $h$ differ,
then $g_m\dots g_1 x$ and $h_n\dots h_1 x$ are clearly different,
while if $m=n$ and these sequences are the same, then
there must be some $j\in M$ such that the sequences of elements
of $G_j$ differ; and if we look at the first terms where this
occurs, then by right cancellativity, the corresponding
terms of~(\ref{x.gr1grm}) will differ, and by our choice of $x_{1,j},$
the elements $gx$ and $hx$ will differ at that step.
\end{proof}

(The above action will, in fact, be strongly faithful -- we could
have handled any finite set of elements of
$\coprod_{i\in I} M_i$ as we did $\{g,h\};$ but we
only need the two-element case, and it allowed simpler notation.)

Now letting $G_0,\ G_1$ be two isomorphic copies of the group $\Sym,$
disjoint except for the identity, each represented naturally
on $\Omega,$ we get as in the paragraph before Lemma~\ref{L.fthfl}
strongly faithful actions of $G_0$ and $G_1$ on
$1\cup\Omega\cup\Omega^2\cup\Omega^3\cup\dots,$ which has cardinality
$\cd.$
Calling this set, regarded as a $\!G_0\!$-set $\Omega_0,$
and regarded as a $\!G_1\!$-set $\Omega_1,$ the above
lemma tells us that the coproduct of $G_0$ and $G_1$ acts faithfully
on $\Omega_0\boxtimes\Omega_1,$ which is also of cardinality
$\cd,$ completing our alternative proof of Lemma~\ref{L.double}(i).

More generally, the submonoid of {\em surjective} endomaps of $\Omega$
satisfies~(\ref{x.hg=g}), so the above construction embeds
the coproduct of two copies of that monoid in $\Se.$

On the other hand, if we delete the right cancellativity
assumption~(\ref{x.hg=g}) from the hypothesis of Lemma~\ref{L.fthfl},
there will in general be no choice of $\!M_i\!$-sets $\Omega_i$ making
$\boxtimes_{i\in I} \Omega_i$ faithful.
For if one of the monoids, say $M_1,$ has elements $a\neq b$
and $c$ satisfying $ac=bc,$ and at least one other of the
monoids is nontrivial, say $M_2,$ with a nonidentity element $d,$ then
in $\coprod M_i$ we find that $adc$ and $bdc$ are distinct
elements having the same action on $\boxtimes_{i\in I} \Omega_i$
for any family of $\!M_i\!$-sets $\Omega_i.$

\section{Appendix: Some conditions on complete lattices.}\label{S.=>=>}
We noted in the discussion preceding Question~\ref{Q.solset+}
that the conditions we had proved on chains of solution sets in
$\Sym$ and $\Se$ were consequences of the
stronger statement that the lattice of all such solution sets
embeds in the system of closed sets
of a topological space having a basis of $\leq\cd$ open sets.
The same observation holds for the results obtained in later
sections on solution sets in $\Eq,$ etc..

The next lemma compares these and related conditions.
Here a {\em $\!\kappa\!$-generated}
topological space means a topological space having a basis
of open sets of cardinality $\leq\kappa\<;$ equivalently,
having a {\em subbasis} of open sets of that cardinality;
equivalently, having such a basis or subbasis of {\em closed} sets.
An {\em embedding} of lattices, of complete lower semilattices, etc.,
means a one-to-one homomorphism of such structures; an embedding
of {\em partially ordered sets} means a (necessarily one-to-one)
map that preserves both the relations $\leq$ and~$\not\leq.$

\begin{lemma}\label{L.=>=>}
Let $\kappa$ be an infinite cardinal, and $A$ a complete lattice.
Then of the following conditions, each implies the next, and conditions
with the same roman numeral and different suffixes are equivalent.

\textup{(i.a)}~ $A$ is embeddable as a complete lower semilattice in
$\Pw[\kappa].$

\textup{(i.b)}~ $A$ is generated as a complete upper semilattice by
a set of $\leq\kappa$ elements.

\textup{(ii.a)}~ $A$ is embeddable as a complete lower semilattice in
a complete lower semilattice generated by $\leq\kappa$ elements.

\textup{(ii.b)}~ $A$ is embeddable as a partially ordered set in
a complete lower semilattice generated by $\leq\kappa$ elements.

\textup{(ii.c)}~ $A$ is embeddable as a partially ordered set in
$\Pw[\kappa].$

\textup{(ii.d)}~ $A$ is embeddable as a complete lower semilattice in
the system of closed subsets
of a $\!\kappa\!$-generated topological space.

\textup{(ii.e)}~ $A$ is embeddable as a partially ordered set in
the system of closed subsets
of a $\!\kappa\!$-generated topological space.

\textup{(ii.a*)-(ii.e*)}~ The duals of \textup{(ii.a)-(ii.e)}; i.e.,
the corresponding statements with ``upper semilattice''
replaced by ``lower semilattice'', and ``closed subsets'' by
``open subsets'', wherever applicable.
\textup{(}So, no change in \textup{(ii.c))}.

\textup{(iii)}~ Every chain in $A$ has a dense subset of
cardinality $\leq\kappa.$

\textup{(iv.a)}~ No chain in $A$ has a family of $>\kappa$ disjoint
intervals.

\textup{(iv.b)}~ No chain in $A$ has $>\kappa$ jumps.

\textup{(v)}~ $A$ contains no well-ordered or reverse-well-ordered
chain of cardinality $>\kappa.$

Moreover, if for each family of conditions whose equivalence
is asserted above, we denote the common condition using the
corresponding roman numeral with suffixes dropped,
then the implications \textup{(i)$\implies$(ii)
$\implies$(iii)} are irreversible for all $\kappa,$
but are both reversible if $A$ is restricted to be a chain;
the implication \textup{(iv)$\implies$(v)} is irreversible for
$\kappa=\aleph_0,$ while the reversibility of
\textup{(iii)$\implies$(iv)} for $\kappa=\aleph_0$ is equivalent to
Suslin's Hypothesis, known to be independent of ZFC.

\end{lemma}\begin{proof}
(i.a)$\implies$(i.b): Given a complete lower semilattice embedding
$f: A\rightarrow\Pw[\kappa],$ let us associate
to each $\alpha\in\kappa$ the meet $g(\alpha)\in A$ of all elements
$x\in A$ satisfying $\alpha\in f(x).$
We see that $g(\alpha)$ will be the least $y$ such that
$\alpha\in f(y),$ and we deduce that every $x\in A$ is the join
in $A$ of $\{g(\alpha)\mid\alpha\in f(x)\}.$
So $\{g(\alpha)\mid\alpha\in\kappa\}$ generates $A$ as a complete
upper semilattice.

(i.b)$\implies$(i.a):
Given a generating set $\{g_\alpha\mid\alpha\in\kappa\}$
for $A$ as a complete upper semilattice, we find that
an embedding $A\rightarrow\Pw[\kappa]$ as complete lower semilattices
is given by the map $x\mapsto\{\alpha\mid g_\alpha\leq x\}.$

(i.a)$\implies$(ii.c) is immediate.

To prove the equivalence of the versions of~(ii), we shall show
(ii.a)$\implies$(ii.d)$\implies$(ii.e)%
$\implies$(ii.b)$\implies$(ii.c)$\implies$(ii.a).
Since~(ii.c) is self-dual, it will follow that these conditions
are also equivalent to their starred variants.

(ii.a)$\implies$(ii.d): It suffices to show
that every complete lower semilattice $A'$ generated by a set
$\{x_\alpha\mid\alpha\in\kappa\}$
is embeddable as a complete lower semilattice
in the system of closed subsets of a $\!\kappa\!$-generated topology.
Given such an $A',$ define for each $x\in A'$ the ``principal downset''
\begin{xlist}\item\label{x.D=}
$D(x)\ =\ \{y\in A'\mid y\leq x\}\ \subseteq\ \Pw[A'].$
\end{xlist}
We see that these sets form a complete lower semilattice isomorphic to
$A';$ hence if we define a topology on the underlying set of $A'$ using
the $D(x)$ as a subbasis of closed sets, $A'$ embeds as a complete lower
semilattice in the complete lattice of closed sets of that topology.
Moreover, the closed sets $D(x_\alpha)$ $(\alpha\in\kappa)$
also form a subbasis of closed sets
for this topology, so it is $\!\kappa\!$-generated, as required.

(ii.d)$\implies$(ii.e)$\implies$(ii.b): Trivial.

(ii.b)$\implies$(ii.c): Suppose $A$ is embeddable as partially
ordered set in a complete lower semilattice $A'$ as in (ii.b).
Thus, $A'$ satisfies the {\em dual} of (i.b).
Hence it satisfies the dual of~(i.a), hence $A,$ being embeddable
in $A'$ as a partially ordered set, satisfies~(ii.c).

(ii.c)$\implies$(ii.a):  Note that the map $D$ of~(\ref{x.D=})
(with $A$ in the role of $A')$
is an embedding $A\rightarrow\Pw[A]$ as complete lower semilattices.
Hence it will suffice to show that (ii.c) implies that the image of
$A$ under $D$ lies in a complete lower subsemilattice
of $\Pw[A]$ generated by $\leq\kappa$ elements.
Given an embedding $f:A\rightarrow\Pw[\kappa]$ as
in~(ii.c), let us define for each $\alpha\in\kappa$
the set $c_\alpha=\{x\in A\mid\alpha\notin f(x)\}\in\Pw[A].$
We see that for every $x\in A,$ the $D(x)$
of~(\ref{x.D=}) is the intersection of those members of
$\{c_\alpha\mid\alpha\in\kappa\}\cup\{A\}\subseteq\Pw[A]$
that contain $x,$ so $D$ carries $A$ into the complete
lower subsemilattice of $\Pw[A]$ generated by these elements.

(ii)$\not\Longrightarrow$(i):
Let $A$ be a lattice consisting of a least element~$0,$
a greatest element~$1,$ and $2^\kappa$ pairwise incomparable
elements lying between these.
Clearly it does not satisfy~(i.b).
To see that it satisfies~(ii.d),
recall that the product topology
on $\Pw[\kappa]$ is $\!\kappa\!$-generated, and let us map
$A$ into the closed sets of that topology by sending $0$ to $\emptyset,$
$1$ to the improper subset, and the $2^\kappa$
intermediate elements to the singletons $\{S\}\subseteq\Pw[\kappa]$
$(S\subseteq\kappa).$
(Alternatively, one can show that $A$ satisfies~(ii.c) by
noting that
$\Pw[\kappa\times 2]\cong\Pw[\kappa]$ has an antichain
of cardinality $2^\kappa,$ consisting of the sets
$(s\times\{0\})\cup (s^c\times\{1\})$ $(s\subseteq\kappa).)$

(ii.c)$\implies$(iii):  It suffices to show that every chain $C$
in $\Pw[\kappa]$ has a dense subset of cardinality $\leq\kappa.$
Let $C$ be such a chain, and for each pair
$\alpha,\,\beta\in\kappa$ such that some $x\in C$
contains $\alpha$ but not $\beta,$ choose such an
element, $x_{\alpha,\beta}.$
This gives a family of $\leq\kappa$ elements which is easily
seen to be dense.

(iii)$\not\Longrightarrow$(ii): Take $A$ as in
the example used to show (ii)$\!\not\Longrightarrow\!$(i),
but this time with $>2^\kappa$ pairwise incomparable elements.
This clearly satisfies~(iii), but in view of its cardinality,
cannot satisfy~(ii.c).

On the other hand, the assertion that if $A$ is a chain
the implications (i)$\!\implies\!$(ii)$\!\implies\!$(iii) are
reversible follows from the obvious implication
(iii)$\!\implies\!$(i.b) in this case.

(iii)$\implies$(iv.a): Given any chain $C$ in $A$ and any dense
set $S$ of $\leq\kappa$ elements of $C,$ we see that for every interval
$[x,y]$ in $C,$ at least one element of $S$ must belong
to $[x,y],$ showing that $C$ cannot have $>\kappa$ disjoint intervals.

Concerning the reverse implication,
recall that Suslin's Hypothesis says
that every totally ordered set $S$ having no uncountable family of
disjoint intervals has a countable dense subset, and that
this is independent of ZFC~\cite{SH}.
Assuming Suslin's Hypothesis, we immediately get (iv.a)$\implies$(iii)
for $\kappa=\aleph_0$ by applying that statement to an arbitrary
chain $C$ in $A.$

In proving the converse assertion, note that~(iv.a) and~(iii)
are statements about a {\em complete} lattice $A$ (first
sentence of the lemma); so we need to show that
if $S$ is any counterexample to Suslin's hypothesis, we can obtain
from it a complete lattice $A$ that satisfies~(iv.a) but not~(iii).
Given such an $S,$ let $A$ be its completion as a totally ordered set;
I claim that $A$ inherits the properties making $S$ a counterexample to
Suslin's conjecture, and thus gives the desired example.
Indeed, for any infinite dense subset $D\subseteq A,$ if we take an
element of $S$ between every two distinct elements of $D$ (and
throw in the greatest and/or least element of
$S$ if these exist), we get a dense subset $D'$ of $S$ of the same
cardinality; so since $S$ has {\em no} countable dense subset,
neither does $A.$
Likewise, if some chain $C\subseteq A$ had an uncountable family of
disjoint intervals, then for each of these intervals $[x,y]_C$
we could choose $x',y'\in S$ with $x\leq x'<y'\leq y,$ getting
an uncountable family of disjoint intervals $[x',y']_S$ in $S;$
so chains in $A$ inherit from $S$ the nonexistence of such families.

Since~(iv.a) and~(iv.b) are negative statements,
we will prove their equivalence in contrapositive form:

$\!\neg\!$(iv.a)$\implies$$\!\neg\!$(iv.b):
If a chain $C$ in $A$ has a family of disjoint intervals
$[x_\alpha,y_\alpha]$ where $\alpha$ ranges over some
$\lambda>\kappa,$ in then the subchain
$C'=\{x_\alpha,y_\alpha\mid\alpha\in\lambda\},$
the pairs $x_\alpha<y_\alpha$ will be jumps.

$\!\neg\!$(iv.b)$\implies$$\!\neg\!$(iv.a):
If $C$ has $>\kappa$ jumps, let us
associate to each jump $x<y$ the two-element interval $[x,y].$
There is the slight difficulty that distinct jumps may not yield
disjoint intervals: the upper endpoint of one may equal the lower
endpoint of the other.
However, if we take a family of these intervals maximal for the
property of being pairwise disjoint, is easy to verify that this still
has cardinality $>\kappa,$ giving the desired assertion.

(iv.b)$\implies$(v) is clear.

To show that for $\kappa=\aleph_0,$
(v)$\not\Longrightarrow$(iv.b), let $A$ be $\mathbb{R}\times 2,$
lexicographically ordered.
This has continuum many jumps (a jump $(r,0)<(r,1)$
for each real number $r),$ but has no uncountable well-ordered
or reverse well-ordered subsets.
\end{proof}

Some remarks on the above lemma:

It is, of course, condition~(ii) and its consequences (iii)-(v)
that are directly relevant to the results of preceding sections.
I have included~(i) for perspective.

(i) is the only condition in the lemma that is not equivalent to its
dual.
To show this inequivalence, let $A\subseteq\Pw[\kappa]$
consist of $\emptyset$ and all sets of cardinality $\kappa.$
This clearly satisfies the dual of~(i.a), but we claim it does
not satisfy (nondualized)~(i.b).
Indeed, given a putative generating set
$\{x_\alpha\mid\alpha\in\kappa\}$ for $A$ as a complete upper
semilattice (where we allow repetitions in the indexing in case
this set has cardinality $<\kappa),$ one can construct by transfinite
recursion an element $y\in\Pw[\kappa]$ having $\kappa$ elements,
but missing at least one element from each $x_\alpha\neq\emptyset$
in our family.
We see that $y$ will belong to $A,$ but not to the complete
upper subsemilattice generated by $\{x_\alpha\}.$

We could have written (i.a*) and (i.b*) for the duals
to (i.a) and (i.b), adding to the lemma their mutual equivalence and
the implication (i*)$\!\implies\!$(ii); but so naming those
conditions would have broken the convention
that conditions beginning with the same roman numeral are equivalent.

Conditions (ii.d) and (ii.e) (and hence their duals), which refer to a
$\!\kappa\!$-generated topological space, are equivalent to the formally
stronger conditions referring
to a $\!\kappa\!$-generated {\em Hausdorff}
(and if we wish, {\em totally disconnected}) topological space.
For given a topological space $X$ as in one of those
statements, which we may assume without loss of generality to be
$\mathrm{T}_0,$ and which
has a subbasis of $\leq\kappa$ closed sets, we can throw in the
complements of those sets to get a stronger topology on $X$
which is still $\!\kappa\!$-generated, but is now totally disconnected
and Hausdorff, and whose lower semilattice of closed sets contains the
lower semilattice of sets closed in the original topology.

It is curious that the example we gave for
(ii)$\!\not\Longrightarrow\!$(i) in the proof of the
lemma is an instance of~(\ref{x.meetz}) with $\cd[I]>\kappa.$
Thus, although the lattice $A$ in question satisfies~(ii),
the lattice $L^{\strt=}_{A_\wedge,1}$ does not even satisfy~(v).
The example given above satisfying the dual
of (i.a), but not (i.b), similarly
contains such an instance of~(\ref{x.meetz}), in view of Sierpi\'nski's
result \cite{WS} that there exists a family of $>\kappa$ subsets
of $\kappa,$ each having cardinality $\kappa,$ but with pairwise
intersections all of smaller cardinality.
These observations suggest the first part of the next question.
The second part is also natural, in view of the simpleminded
example we used for~(iii)$\!\not\Longrightarrow\!$(ii).

\begin{question}\label{Q.L|-v}
\textup{(a)} If we add to the hypotheses of Lemma~\ref{L.=>=>}
the assumption that $L^{\strt=}_{A_\wedge,1}$ satisfies \textup{(v)},
does this change the validity of the nonimplications shown?
\textup{(}If not, we might try imposing the stronger
condition \textup{(iv)} or \textup{(iii)} on these
lattices, and/or looking at $L^{\strt=}_{A_\wedge,\,J}$ for $J$ of
larger cardinality, up to $\kappa.)$

\textup{(b)} If we add to the hypotheses of Lemma~\ref{L.=>=>}
the assumption that $A$ has cardinality $\leq 2^\kappa,$
or that it has order-dimension $\leq\kappa$ \textup{(}i.e.,
is embeddable as a partially ordered set in a direct
product of $\leq\kappa$ totally ordered sets -- both
conditions being implied by \textup{(ii))}, does this affect the
validity of the assertion \textup{(iii)$\not\Longrightarrow$(ii)}?
\end{question}

Our final corollary, below, answers a couple of other questions
suggested by that lemma.
The formulation of statement~(b.2) of that corollary uses implicitly
the fact that
in a complete upper semilattice $A'$ with least element, every subset
has a greatest lower bound, so that $A'$ may be regarded as
a complete lattice (though if $A'$ was obtained as a complete upper
subsemilattice of a complete lattice $A,$ the meet operation
of $A',$
will not in general agree with that of $A).$
Statement~(b.2*) uses the dual observation.

\begin{corollary}\label{C.=>=>}
\textup{(a)}~ For the conditions of Lemma~\ref{L.=>=>} that treat
$A$ only as a partially ordered set, namely
\textup{(ii.b), (ii.c), (ii.e), (ii.b*), (ii.c*), (ii.e*)}
and \textup{(iii)-(v)}, the implications stated in that lemma
\textup{(}for $A$ the underlying partially ordered
set of a complete lattice\textup{)}
in fact hold for any partially ordered set $A.$

\textup{(b)}~ Let $C$ be a complete lattice, and $A$ a nonempty
subset of $C.$
Then the following conditions are equivalent.\\
\textup{(b.1)}~ As a partially ordered set, $A$ forms
a complete lattice.\\
\textup{(b.2)}~ $A$ has a greatest element, and forms a complete lower
subsemilattice of a complete upper subsemilattice
$B$ of $C$ having a least element.\\
\textup{(b.2*)}~ $A$ has a least element, and forms a complete upper
subsemilattice of a complete lower subsemilattice
$B$ of $C$ having a greatest element.

\textup{(}Note: in \textup{(b.2)} and \textup{(b.2*)}, the statement
that $A$ is a complete lower or upper subsemilattice
of $B$ means that it is closed under the relevant
meet or join operations of $B,$ not under those of $C,$ which do not
in general carry $B$ into itself.\textup{)}
\end{corollary}\noindent
{\em Sketch of proof.}
(a):  The versions of condition (ii) listed here all involve embeddings
of $A$ as a partially ordered set in a certain complete semilattice
$B,$ which, by adjoining a greatest or least element if
necessary, can be assumed a complete lattice.
The mutual equivalence of these conditions for complete lattices,
applied to $B,$ gives embeddings of $B;$ these in turn
lead to the embeddings of $A$ as a partially ordered set that we want.
The same method yields the implication (ii)$\implies$(iii).
The proofs of the remaining implications were entirely order-theoretic,
and go over unchanged.

(b):  Clearly~(b.2)$\!\implies\!$(b.1).
On the other hand, if $A$ satisfies~(b.1), and
one lets $B$ denote the result of closing $A$ in $C$
under arbitrary joins and throwing in the
least element, then one finds that arbitrary meets in $A$
are still meets under the operation of $B,$ so that $B$ witnesses~(b.2).
Thus, (b.1)$\!\iff\!$(b.2), and by duality,~(b.1)$\!\iff\!$(b.2*).\qed

George M. Bergman\\
Department of Mathematics\\
University of California\\
Berkeley, CA 94720-3840\\
USA\\
gbergman@math.berkeley.edu
\end{document}